\documentclass[preprint,11pt]{elsarticle}

\usepackage[a4paper, total={6in, 8in}]{geometry}

\makeatletter
\def\ps@pprintTitle{%
 \let\@oddhead\@empty
 \let\@evenhead\@empty
 \def\@oddfoot{\centerline{\thepage}}%
 \let\@evenfoot\@oddfoot}
\makeatother

\usepackage{amssymb}
\usepackage{amsthm}
\usepackage{amsmath}

\usepackage{adjustbox}
\usepackage{empheq}

\usepackage{tikz,xcolor}
\usepackage{pgfplots}
\pgfplotsset{compat=newest}
\usepackage{mathtools}
\usepackage{mathrsfs}
\usepackage{enumitem}

\usepackage{graphicx}
\usepackage{subcaption}

\usetikzlibrary{patterns}

\newtheorem{remark}{Remark}
\newtheorem{theorem}{Theorem}

\newtheorem{lemma}{Lemma}
\newtheorem{proposition}{Propositon}
\newtheorem{assumption}{Assumption}

 \usepackage{pgfplotstable} 
 \usepackage{booktabs}      

\usepackage{stmaryrd}
\usepackage{tikz}
\usetikzlibrary{positioning}
\usetikzlibrary{shapes,arrows}

\biboptions{sort&compress}

\usepackage{listings}
\usepackage{color} 
\usepackage[normalem]{ulem} 

\newcommand{\R}{\mathbb{R}}

\makeatletter
\def\smallunderbrace#1{\mathop{\vtop{\m@th\ialign{##\crcr
   $\hfil\displaystyle{#1}\hfil$\crcr
   \noalign{\kern3\p@\nointerlineskip}%
   \tiny\upbracefill\crcr\noalign{\kern3\p@}}}}\limits}
\makeatother

\pgfmathsetmacro{\width} {6cm}
\pgfmathsetmacro{\height} {6cm}

\begin{document}

\begin{frontmatter}

\title{Multistage DPG time-marching scheme for nonlinear problems}

\author[1,2]{Judit Mu\~noz-Matute}
\author[2]{Leszek Demkowicz}

\address[1]{Basque Center for Applied Mathematics (BCAM), Bilbao, Spain}
\address[2]{Oden Institute for Computational Engineering and Sciences,\\ The University of Texas at Austin, Austin, USA}

\begin{abstract}
In this article, we employ the construction of the time-marching Discontinuous Petrov-Galerkin (DPG) scheme we developed for linear problems to derive high-order multistage DPG methods for non-linear systems of ordinary differential equations. The methodology extends to abstract evolution equations in Banach spaces, including a class of nonlinear partial differential equations. We present three nested multistage methods: the hybrid Euler method and the two- and three-stage DPG methods. We employ a linearization of the problem as in exponential Rosenbrock methods, so we need to compute exponential actions of the Jacobian that change from time step to time step. The key point of our construction is that one of the stages can be post-processed from another without an extra exponential step. Therefore, the class of methods we introduce is computationally cheaper than the classical exponential Rosenbrock methods. We provide a full convergence proof to show that the methods are second, third, and fourth-order accurate, respectively. We test the convergence in time of our methods on a 2D + time semi-linear partial differential equation after a semidiscretization in space. 
\end{abstract}


\end{frontmatter}

\section{Introduction}\label{Sec:Intro}

Exponential time integrators \cite{hochbruck2010exponential} are a class of time-marching schemes for solving stiff systems of non-linear or semi-linear Ordinary Differential Equations (ODEs), usually arising from the spatial discretization of Partial Differential Equations (PDEs). These methods are of particular interest as they allow the use of large time steps because of their capacity to handle the stiffness of the problem. There exist different types of exponential integrators, including exponential Runge-Kutta methods \cite{hochbruck2005explicit,hochbruck2005exponential,luan2014explicit}, exponential Rosenbrock methods \cite{hochbruck2009exponential,caliari2009implementation,luan2014exponential}, exponential multistep methods \cite{hochbruck2011exponential,calvo2006class}, or Lawson methods \cite{ostermann2020lawson,krogstad2005generalized}, among others. All these integrators require to compute the action of exponential-related functions (called $\varphi$-functions) over vectors. Although the first articles about exponential integrators can be traced back to the 70's \cite{friedli2006verallgemeinerte}, they have attracted more attention in the last decade due to the recent developments in numerical linear algebra on efficient algorithms to compute these actions \cite{berland2007expint,higham2020catalogue,al2011computing,higham2005scaling,niesen2012algorithm}, making these types of methods significantly more competitive. 

On the other hand, the Discontinuous Petrov-Galerkin (DPG) method with optimal test functions \cite{demkowicz2014overview,demkowicz2010class,demkowicz2011class} is a well-established method within the class of stabilized Finite Element Methods \cite{zienkiewicz2005finite} to approximate the solution of challenging PDEs. This approach has presented an impressive performance in approximating challenging problems in engineering like advection-dominated diffusion \cite{demkowicz2013robust}, Stokes flows \cite{roberts2014dpg}, linear elasticity \cite{bramwell2012locking}, or time-dependent problems \cite{demkowicz2017spacetime}, among many others. The main idea behind the DPG method is to construct optimal test functions such that the discrete stability of the method is guaranteed. 

In our previous work \cite{munoz2021adpg,munoz2021equivalence,munoz2022error,munoz2022combining}, we applied the DPG theory to stiff systems of linear ODEs arising from the semidiscretization of linear transient PDEs, and we derived a DPG-based time-marching scheme. We proved that if we employ an ultraweak variational formulation in the time variable, the optimal test functions that ensure the stability of the scheme are exponential-related functions that can be expressed in terms of the $\varphi$-functions. In our construction, we have two types of variables in time: traces (point values) and fields (piecewise polynomials). We proved that the DPG time-marching scheme is equivalent to classical exponential integrators for the trace variables, and additionally, it delivers the $L^2$-projection of the exact solution in the interior of the time intervals. In summary, the DPG time-marching scheme can be viewed as the hybridization of classical exponential integrators for linear problems. Apart from the guaranteed stability of the resulting method, the are other benefits when considering a variational formulation in the time variable. One of them is that the DPG method provides a built-in error representation function that can be employed for adaptivity \cite{munoz2022error}. Finally, in goal-oriented adaptive strategies, a variational formulation in time is essential in order to express the error of a pre-defined quantity of interest as an integral over the whole space-time domain \cite{bangerth2010adaptive}.

In this article, we employ the structure of the DPG time-marching scheme we introduced for linear problems to construct a new class of higher-order methods for non-linear problems. We derive three methods: the hybrid Euler method and the two- and three-stage DPG. These methods are of orders two, three, and four, respectively. For that, we first consider a linearization of the problem at the initial time-step $t_n$ of each time interval as in Rosebrock-type methods \cite{hochbruck2009exponential}. We realize that if we approximate the nonlinear remainder with the known value of the solution at $t_n$, we are in the linear setting. Therefore, we can apply the lowest order (piecewise constants) DPG-time marching scheme from \cite{munoz2021equivalence}, and we prove that it is the hybridization of the classical exponential Euler method. Here, we compute the interior variable with an exponential step, and then we post-process the trace variable without computing an extra action of any $\varphi$-function. The hybrid exponential Euler method is second-order and has the same cost as the classical exponential Euler method. Then, we construct methods up to order four employing a hybrid Euler step for defining the internal stages. We define the update for the two-stage DPG method as in classical exponential Runge-Kutta methods, and we derive the so-called stiff order conditions for the coefficient functions to obtain a third-order method. However, for the three-stage DPG method, it is not enough to derive the stiff order conditions in the classical way to obtain a fourth-order method. In addition, we need to add a term to the update that depends upon the Jacobian to obtain the desired order. The key point of the construction of the fourth-order method is that one of the stages is computed as a postprocessing of the other one without an extra exponential action, being therefore cheaper than the classical exponential Rosenbrock methods. Moreover, the third- and fourth order methods are nested so they can be employed to adapt the time-step size. 

The article is organized as follows: Section \ref{Sec:Model} describes the model problem, the linearization we consider, and the assumptions of the solution and the non-linear term. Section \ref{Sec:Linear} summarizes the idea of applying the DPG method as a time-integrator for linear problems. Section \ref{Sec:Methods} introduces the construction of the hybrid Euler method and the two and three-stage DPG methods. In Section \ref{Sec:Stiff}, we derive the stiff order conditions by analyzing the local truncation error, and in Section \ref{Sec:Convergence}, we provide the final convergence proofs. In Section \ref{Sec:Results}, we test the convergence of the methods on a benchmark 2D+time semilinear problem. We also numerically test some energy preservation properties of the methods in 1D+time semilinear and quasilinear problems. Section \ref{Sec:Conclusions} discusses the conclusions and future work. Finally, in \ref{App:Summarize} we provide a summary of the notation we employed throughout the article. 

\section{Model problem, assumptions, and linearization}\label{Sec:Model}
Letting $I=[0,T]\subset\R$, we consider the following first-order nonlinear autonomous system of Ordinary Differential Equations (ODEs)
\begin{equation}\label{ODE1}
\displaystyle{ \left\{
\begin{split}
u'(t)&=F(u(t)),\\
u(0)&=u_{0}.\\
\end{split}
\right.} 
\end{equation}
We consider the framework of strongly continuous semigroups \cite{pazy2012semigroups} for the convergence analysis. Therefore, letting $X$ be a Banach space with norm $\Vert\cdot\Vert$, we consider the following usual assumptions \cite{hochbruck2009exponential} on the solution $u:[0,T)\longrightarrow X$ and the nonlinear term in (\ref{ODE1}).
\begin{assumption}\label{Ass1}
We assume that $u\in C([0,T);X)\cap C^1((0,T);X)$
and that $F:X\longrightarrow X$ is sufficiently many times Fr\'echet differentiable in a strip along the exact solution with all derivatives uniformly bounded. Therefore, the Jacobian operator $J(u)=\frac{\partial F}{\partial u}(u)$ satisfies the Lipschitz condition 
\begin{equation}\label{Ass1for}
\Vert J(u)-J(v)\Vert_{_{\mathcal{L}(X)}}\leq C\Vert u-v\Vert,
\end{equation}
in a neighborhood of the exact solution $u$.
\end{assumption}
\begin{assumption}\label{Ass2}
The linear operator $J$ is the generator of a strongly continuous semigroup\footnote{Note that this assumption limits the class of possible nonlinearities we can address with this framework. We refer to \cite{hochbruck2005explicit}, where the authors construct explicit exponential Runge-Kutta methods for more general nonlinearities by assuming that the linear operator generates an analytic semigroup.} on $X$, i.e., there exist constants $C$ and $\omega$ such that 
\begin{equation}\label{Ass2for}
\lVert e^{tJ}\rVert_{_{\mathcal{L}(X)}}\leq Ce^{\omega t},\;\;\;t\geq0,
\end{equation}
holds uniformly in a neighborhood of the exact solution of (\ref{ODE1}).
\end{assumption}

For the construction of the methods, we first consider a uniform partition of the time interval
\begin{equation}\label{TimePartition}
0=t_{0}<t_{1}<\ldots<t_{N-1}<t_{N}=T,
\end{equation}
we define $I_{n}=(t_{n},t_{n+1})$ and $h=t_{n+1}-t_{n},\;\forall n=0,\ldots,N-1$. Letting $u_n$ be the numerical approximation to the solution of (\ref{ODE1}) at $t_n$, i.e. $u_n\approx u(t_n)$, we rewrite (\ref{ODE1}) as a semilinear system as follows
\begin{equation}\label{linODE1}
\displaystyle{ \left\{
\begin{split}
u'(t)&=J_nu(t)+g_n(u(t)),\\
u(0)&=u_{0},\\
\end{split}
\right.} 
\end{equation}
where 
\begin{equation}\label{JacRem}
J_n=\frac{\partial F}{\partial u}(u_n),\;\;\;g_n(u)=F(u)-J_nu,
\end{equation}
are the Jacobian of $F$ evaluated at $u_n$ and the nonlinear remainder, respectively. 

\begin{remark}\label{rmk:1}
If $F$ is already semilinear, i.e., $F(u)=Au+f(u)$, with $A$ being a linear operator, we have $J_n=A+\frac{\partial f}{\partial u}(u_n)$ and $g_n(u)=f(u)-\frac{\partial f}{\partial u}(u_n)u$. Note that the framework presented here also covers the case where $A$ is a differential operator, and $f$ is an algebraic function of $u$. Therefore, the methods developed in this article can be applied to a wide range of Partial Differential Equations (PDEs). In the numerical results in Section \ref{Sec:Results}, we consider this case being $A$ a discretization of the Laplacian operator.
\end{remark}

\begin{remark}\label{rmk:2}
If the system is non-autonomous, i.e., $u'(t)=F(t,u(t))$, we can easily reduce it to be autonomous with the following change of variables
$$U=\begin{bmatrix}t\\u\end{bmatrix},\;\;\mathcal{F}(U)=\begin{bmatrix}1\\F(t,u)\end{bmatrix}.$$
The linearized system then becomes $U'=\mathcal{J}_nU+\mathcal{G}_n(U)$, where the Jacobian is defined as 
$$\mathcal{J}_n=\begin{bmatrix}0&0\\\frac{\partial F}{\partial t}(t_n,u_n)&\frac{\partial F}{\partial u}(t_n,u_n)\end{bmatrix},$$
and $\mathcal{G}_n(U)=\mathcal{F}(U)-\mathcal{J}_nU$ is the nonlinear remainder. 
\end{remark}

\begin{remark}
Note that
\begin{equation}\label{vanishgn}
\frac{\partial g_n}{\partial u}(u_n)=\frac{\partial F}{\partial u}(u_n)-J_n=0.
\end{equation}
As we will see in Section \ref{Sec:Stiff}, this property implies that all the approximation methods derived for (\ref{linODE1}) are at least second order by construction. 
\end{remark}

The integral representation of the solution of (\ref{linODE1}) at $t_{n+1}$, also known as the \textit{variation-of-constants} formula, reads
\begin{equation}\label{VOC}
u(t_{n+1})=e^{hJ_n}u(t_n)+h\int_0^1e^{(1-\theta)hJ_n}g_n(u(t_n+h\theta))d\theta,
\end{equation}
and different approximations of the nonlinear term in (\ref{VOC}) lead to different exponential time integration methods. All these methods are expressed in terms of the so-called \textit{$\varphi$-functions} defined as
\begin{equation}\label{PhiFunctions}
\displaystyle{\left\{
\begin{split}
\varphi_{0}(z)&=e^{z},\\
\varphi_{p}(z)&=\int_{0}^{1}e^{(1-\theta)z}\frac{\theta^{p-1}}{(p-1)!}d\theta,\;\forall p\geq1,
\end{split}
\right.} 
\end{equation}
which satisfy the following recurrence relation
\begin{equation}\label{Recurrence}
\varphi_{p+1}(z)=\frac{1}{z}\left(\varphi_{p}(z)-\frac{1}{p!}\right).
\end{equation}

\section{Overview of the DPG time-marching scheme for linear problems}\label{Sec:Linear}
In our previous works on the DPG method for linear transient problems \cite{munoz2021adpg,munoz2021equivalence,munoz2022error,munoz2022combining}, we proved that from an ultraweak variational formulation in time, we obtain a decoupled time-marching scheme where: (a) we obtain the variation-of-constants formula for the solution at each time instant $t_n$ (traces), and (b) we additionally compute the solution at each time interval $I_n$ (fields). We will employ this hybrid approximation in time to construct the internal stages of a Runge-Kutta-like method in the next section. To be self-contained, we present here a brief overview of the DPG-based time-marching scheme with the basic concepts we will employ in this article.

\subsection{DPG time-marching scheme for linear systems}
\begin{figure}[h]
\centering
\begin{tikzpicture}
\draw[thick][-] (0,0) -- (3,0) ; 
		\foreach \x in {0,3}
     		\draw[thick][-]  (\x,-0.1) -- (\x,0.1) node[anchor=north] {};
		\draw (0,-0.35) node{\small$t_n$};
		\draw (3,-0.35) node{\small$t_{n+1}$};

		\draw [thick] (0,0.8) .. controls (1,0.2) and (2,1.8) ..  (3,1.25) ;
		
		\draw plot[mark=*, mark options={fill=black}, mark size=1pt] (0,0.5);
		\draw plot[mark=*, mark options={fill=white}, mark size=1pt] (0,0.8);
		\draw plot[mark=*, mark options={fill=white}, mark size=1pt] (3,1.25);
		\draw plot[mark=*, mark options={fill=black}, mark size=1pt] (3,0.7);
		
		\draw (0-0.35,0.5) node{\small$u_{n}$};
		\draw (1.5,1.5+0.2) node{\small$\tilde{u}_{n}(t)$};
		\draw (3+0.6,0.7) node{\small$u_{n+1}$};
\end{tikzpicture}
\caption{Hybrid approximation of problem (\ref{ODE2}) at $I_n$.}
\label{SolutionPG}
\end{figure}
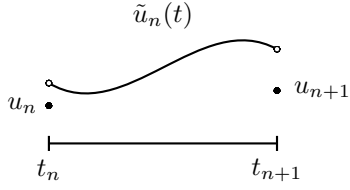

We consider the following linear system of ODEs 
\begin{equation}\label{ODE2}
\displaystyle{ \left\{
\begin{split}
u'(t)+Au(t)&=f(t),\\
u(0)&=u_{0}.\\
\end{split}
\right.} 
\end{equation}
In \cite{munoz2021equivalence}, we considered an ultraweak variational formulation for problem (\ref{ODE2}). For that, we multiply the equation by suitable test functions, integrate by parts, and introduce additional variables for the traces in time (hybridization). We also proved that the optimal test functions to consider are the ones satisfying the adjoint equation (exponential-related functions). Finally, for each time interval $I_n$, we obtain an approximation $\{\tilde{u}_n(t),u_{n+1}\}$ like in Figure \ref{SolutionPG}, where $\tilde{u}_n(t)=\displaystyle{\sum_{j=0}^p\tilde{u}^j_n\left(\frac{t-t_n}{h}\right)^j}$ is a polynomial of order $p$ and $u_{n+1}$ is the point-value (trace) approximation at $t_{n+1}$. 

Finally, we obtained the following decoupled time-marching scheme to compute $u_{n+1}$ and the coefficients of $\tilde{u}_n(t)$
\begin{equation}\label{DPGpSys}
\displaystyle{\left\{
\begin{split}
\sum_{j=0}^{p}\frac{h}{r+j+1}\tilde{u}_n^j&=\tilde{v}_{n}^r(A,t_{n}) u_{n}+\int_{I_n}\tilde{v}^r_{n}(A,t)f(t)dt,\;\forall r=0,\ldots,p,\\
u_{n+1}&=v_n(A,t_{n})u_{n}+\int_{I_{n}}v_n(A,t)f(t)dt,\\
\end{split}
\right.} 
\end{equation}
where $\{v_n,\tilde{v}_n^r,\;r=0,\ldots,p\}$ are the optimal test function. Here, the functions $\tilde{v}_{n}^r(z,t)$ and $v_n(z,t)=e^{z(t-t_{n+1})}$, satisfy the following recurrence formula
\begin{equation}\label{OptTest}
\begin{split}
\tilde{v}_{n}^r(z,t)&=\frac{1}{z}\left(\left(\frac{t-t_{n}}{h}\right)^{r}+\frac{r}{h}\tilde{v}^{r-1}_{n}(z,t)-v_n(z,t)\right),\;\forall r=0,\ldots,p.
\end{split}
\end{equation}
In (\ref{DPGpSys}), we obtain the variation-of-constants formula in the second equation and a scheme that delivers the $L^2-$projection of the analytical solution in the element interior in time. 

\subsection{Post-processing of the trace variables}
Note that from (\ref{OptTest}) with $r=0$, we have
$$v_n(z,t)=1-z\tilde{v}_{n}^0(z,t),$$
and we can express system (\ref{DPGpSys}) as
\begin{equation}\label{DPGpSysPost}
\displaystyle{\left\{
\begin{split}
\sum_{j=0}^{p}\frac{h}{r+j+1}\tilde{u}_n^j&=\tilde{v}_{n}^r(A,t_{n})u_{n}+\int_{I_n}\tilde{v}^r_{n}(A,t)f(t)dt,\;\forall r=0,\ldots,p,\\
u_{n+1}&=u_{n}-\sum_{j=0}^p\frac{h}{j+1}A\tilde{u}_n^j+\int_{I_{n}}f(t)dt.\\
\end{split}
\right.} 
\end{equation}

Therefore, we can first compute the field variables and then we can post-process the trace variable without an extra action of any exponential-reated function \footnote{For example, for $p=0$ and approximating the source with a constant value $f(t)\approx f(t_n)$, we compute the field variable with one exponential step: $\tilde{u}_n^0=\varphi_1(-hA)u_n+h\varphi_2(-hA)f(t_n)$. Then, we post-process for the trace variable with the equation $u_{n+1}=u_n-hA\tilde{u}_n^0+hf(t_n)$, without an extra action of a $\varphi-$function.}. The resulting scheme (\ref{DPGpSysPost}) delivers an approximation of the solution that is equivalent to classical exponential integrators in the trace variables plus an approximation in the interiors. In the next section, we will employ \eqref{DPGpSysPost} with $p=0$ as the internal stages to construct our fourth-order multistage method. Note that the key idea is that the second equation in \eqref{DPGpSysPost} does not require the integration of any exponential-related function.

\subsection{Relation between the optimal test functions and the $\varphi$-functions}
In order to implement scheme (\ref{DPGpSysPost}), we showed in \cite{munoz2021equivalence} the relation between the optimal test functions and the $\varphi-$functions. Given the optimal test functions over the master element (0,1), $\{v,\tilde{v}^r,\;r=0,\ldots,p\}$, we have 
\begin{equation}\label{ToMaster}
\begin{split}
v_n(z,t_{n}+\theta h)&=v(zh,\theta),\\
\tilde{v}^{r}_{n}(z,t_{n}+\theta h)&=h\tilde{v}^r(zh,\theta),\;\forall r=0,\ldots,p,
\end{split}
\end{equation}
and also 
\begin{equation}\label{RelationOptPhi}
\begin{split}
\tilde{v}^{r}(z,0)&=\sum_{j=0}^{r}\frac{r!}{j!}(-1)^{r-j}\varphi_{r-j+1}(-z),\\
\int_{0}^{1}\tilde{v}^{r}(z,t)\frac{t^{q}}{q!}dt&=\sum_{j=0}^{r}\frac{r!}{j!}(-1)^{r-j}\varphi_{r-j+q+2}(-z).
\end{split}
\end{equation}
We can now easily integrate (\ref{DPGpSysPost}) over the master element and then employ (\ref{ToMaster}) and (\ref{RelationOptPhi}) to express the DPG time-marching scheme in terms of $\varphi$-functions. In general, if we select $s$ integration points to approximate the source term in (\ref{DPGpSysPost}), we set $p=s-1$ to obtain optimal convergence rates. 

\section{Multistage DPG time-marching scheme for nonlinear problems}\label{Sec:Methods}
In this section, we introduce several multistage Runge-Kutta-like methods (up to order $4$) for approximating the semilinear problem (\ref{linODE1}) employing the lowest-order ($p=0$) DPG construction for the internal stages.  

\subsection{Hybrid exponential Euler method}
The simplest choice is to approximate the nonlinear remainder in (\ref{linODE1}) and (\ref{VOC}) by its value at $u_n$, which is known from the previous time step, i.e., 
$$g_n(u(t))\approx g_n(u_n),\;\forall t\in I_n,$$
and we are within the framework of linear problems introduced in Section \ref{Sec:Linear}. Therefore, we can employ the lowest order DPG method ($p=0$) and we obtain from (\ref{DPGpSysPost}) and properties (\ref{Recurrence}) and (\ref{RelationOptPhi}) the following scheme 
\begin{equation}\label{ExpEuler}
\displaystyle{\left\{
\begin{split}
\tilde{u}^0_{n}&=u_{n}+h\varphi_2(hJ_n)(g_n(u_{n})+J_nu_{n}),\\
u_{n+1}&=u_{n}+hJ_n\tilde{u}^0_n+hg_n(u_{n}).
\end{split}
\right.} 
\end{equation}

\begin{figure}[h]
\centering
\begin{tikzpicture}
\draw[thick][-] (0,0) -- (3,0) ; 
		\foreach \x in {0,3}
     		\draw[thick][-]  (\x,-0.1) -- (\x,0.1) node[anchor=north] {};
		\draw (0,-0.35) node{\small$t_{n}$};
		\draw (3,-0.35) node{\small$t_{n+1}$};
		
		\draw [thick] (0,1.25) -- (3,1.25) ;
		
		\draw plot[mark=*, mark options={fill=black}, mark size=1pt] (0,0.5);
		\draw plot[mark=*, mark options={fill=white}, mark size=1pt] (0,1.25);
		\draw plot[mark=*, mark options={fill=white}, mark size=1pt] (3,1.25);
		\draw plot[mark=*, mark options={fill=black}, mark size=1pt] (3,0.7);
		
		\draw (0-0.1,0.75) node{\small$u_{n}$};
		\draw (1.5,1.5+0.2) node{\small$\tilde{u}_n^0$};
		\draw (3+0.35,0.9) node{\small$u_{n+1}$};
\end{tikzpicture}
\caption{Solution of the hybrid exponential Euler method (\ref{ExpEuler}).}
\label{SolExpEuler}
\end{figure}
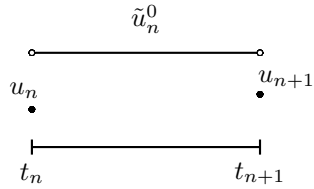

Method (\ref{ExpEuler}) delivers the solution for the trace variable $u_{n+1}$ that coincides with the classical exponential Euler method and, additionally, a constant approximation of the solution $\tilde{u}_n^0$ in the interior of the time interval (see Figure \ref{SolExpEuler}). Here, we only perform a single action of a $\varphi-$function at each time interval. Therefore, with this DPG method, we compute a hybrid approximation in the whole time interval with the same computational cost as for the classical exponential Euler method.

We will show in the next section that method (\ref{ExpEuler}) is second order accurate due to property (\ref{vanishgn}). In the following subsections, we employ approximation (\ref{ExpEuler}) for the internal stages to construct higher-order multistage methods. 

\subsection{Two-stage DPG method}
To construct a two-stage method, we employ the constant approximation we introduced in the first equation of (\ref{ExpEuler}) (denoted here by $u_{n_2}$). Then, we approximate the solution at the next time step employing a linear combination of the nonlinear remainder evaluated at both $u_n$ and $u_{n_2}$ as follows
\begin{equation}\label{generalDPG2}
\displaystyle{\left\{
\begin{split}
u_{n+1}&=e^{hJ_n}u_n+hb_1(hJ_n)g_n(u_{n})+hb_2(hJ_n)g_n(u_{n_2}),\\
u_{n_2}&=u_n+h\varphi_2(hJ_n)(g_n(u_{n})+J_nu_n).\\
\end{split}
\right.} 
\end{equation}
In the next section, we will derive the so-called stiff order conditions for the coefficient functions $b_1$ and $b_2$ to obtain a third-order method. Similar to the classical exponential Runge-Kutta methods, the coefficients $b_1$ and $b_2$ will be expressed in terms of the $\varphi-$functions.

\subsection{Three-stage DPG method}
For the three-stage method, we employ the whole hybrid approximation (\ref{ExpEuler}) as internal stages (denoted by $u_{n_{2}}$ and $u_{n_{3}}$). In this case, however, it is not possible to build a fourth-order method employing solely a linear combination of evaluations of $g_n$. We need to introduce a ``correction factor'' denoted by $C_{n_2}$ below that will be given in terms of the Fr\'echet derivative of $g_n$ in the next section
\begin{equation}\label{generalDPG3}
\displaystyle{\left\{
\begin{split}
u_{n+1}&=e^{hJ_n}u_n+hb_1(hJ_n)g_n(u_{n})+hb_2(hJ_n)(g_n(u_{n_2})+C_{n_2})+hb_3(hJ_n)g_n(u_{n_3}),\\
u_{n_2}&=u_n+h\varphi_2(hJ_n)(g_n(u_{n})+J_nu_n),\\
u_{n_3}&=u_n+hJ_nu_{n_2}+hg_n(u_n).\\
\end{split}
\right.} 
\end{equation}

\section{Stiff order conditions and local truncation error}\label{Sec:Stiff}
In order to derive the stiff order conditions for the two and three-stage DPG methods, we employ a Taylor expansion in the variations-of-constant formula (\ref{VOC}) and compare it against the Taylor expansions for expressions (\ref{generalDPG2}) and (\ref{generalDPG3}).

\subsection{Taylor expansion of the analytical solution}
Letting $\hat{u}_n=u(t_n)$ be the exact solution of $u$ at $t_n$, the linearization of (\ref{ODE1}) at $\hat{u}_n$ now reads 
$$u'(t)=\hat{J}_nu(t)+\hat{g}_n(u(t))$$
where 
$$\hat{J}_n=\frac{\partial F}{\partial u}(\hat{u}_n),\;\;\hat{g}_n(u)=F(u)-\hat{J}_nu.$$

The Taylor expansion of $\hat{g}_n$ at $\hat{u}_n$ is
\begin{equation}\label{Taylorf}
\hat{g}_n(u(t_n+h\theta))=\hat{g}_n(\hat{u}_n)+\frac{1}{2!}\hat{g}''_n(\hat{u}_n)(w_n,w_n)+\frac{1}{3!}\hat{g}'''_n(\hat{u}_n)(w_n,w_n,w_n)+\mathcal{R},
\end{equation}
where $w_n:=u(t_n+h\theta)-\hat{u}_n$ and we used from (\ref{vanishgn}) that $\hat{g}'_n(\hat{u}_n)=0$. Equivalently, the Taylor expansion of $u$ at $t_n$ is
\begin{equation}\label{Tayloru}
u(t_n+h\theta)=\hat{u}_n+\hat{u}'_nh\theta+\frac{1}{2!}\hat{u}''_n(h\theta)^2+\frac{1}{3!}\hat{u}'''_n(h\theta)^3+\mathcal{R}.
\end{equation}

Substituting (\ref{Tayloru}) into (\ref{Taylorf}), we obtain 
\begin{equation}\label{Taylorfu}
\begin{split}
\hat{g}_n(u(t_n+h\theta))&=\hat{g}_n(\hat{u}_n)+\frac{(h\theta)^2}{2!}\hat{g}''_n(\hat{u}_n)(\hat{u}'_n,\hat{u}'_n)+\frac{(h\theta)^3}{2!}\hat{g}''_n(\hat{u}_n)(\hat{u}'_n,\hat{u}''_n)\\
&+\frac{(h\theta)^3}{3!}\hat{g}'''_n(\hat{u}_n)(\hat{u}'_n,\hat{u}'_n,\hat{u}'_n)+\mathcal{R},
\end{split}
\end{equation}
then substituting (\ref{Taylorfu}) into (\ref{VOC})
\begin{equation}\label{Taylorfu2}
\begin{split}
\hat{u}_{n+1}=e^{h\hat{J}_n}\hat{u}_n&+h\int_0^1e^{(1-\theta)h\hat{J}_n}\hat{g}_n(\hat{u}_n)d\theta+h^3\int_0^1e^{(1-\theta)h\hat{J}_n}g_n''(\hat{u}_n)(\hat{u}'_n,\hat{u}'_n)\frac{\theta^2}{2!} d\theta\\
&+h^4\int_0^1e^{(1-\theta)h\hat{J}_n}\left(3\hat{g}_n''(\hat{u}_n)(\hat{u}'_n,\hat{u}''_n)+\hat{g}'''_n(\hat{u}_n)(\hat{u}'_n,\hat{u}'_n,\hat{u}'_n)\right)\frac{\theta^3}{3!} d\theta+\mathcal{R}.\\
\end{split}
\end{equation}

From the definition of the $\varphi$-functions (\ref{PhiFunctions}), we have 
\begin{equation}\nonumber
\begin{split}
\hat{u}_{n+1}=e^{h\hat{J}_n}\hat{u}_n&+h\varphi_1(h\hat{J}_n)\hat{g}_n(\hat{u}_n)+h^3\varphi_3(h\hat{J}_n)\hat{g}''_n(\hat{u}_n)(\hat{u}'_n,\hat{u}'_n)\\
&+h^4\varphi_4(h\hat{J}_n)\left(3\hat{g}_n''(\hat{u}_n)(\hat{u}'_n,\hat{u}''_n)+\hat{g}'''_n(\hat{u}_n)(\hat{u}'_n,\hat{u}'_n,\hat{u}'_n)\right)+\mathcal{O}(h^5),
\end{split}
\end{equation}
and equivalently employing recurrence formula (\ref{Recurrence})
\begin{equation}\label{TaylorExact}
\begin{split}
\hat{u}_{n+1}=\hat{u}_n&+h\varphi_1(h\hat{J}_n)(\hat{g}_n(\hat{u}_n)+\hat{J}_n\hat{u}_n)+h^3\varphi_3(h\hat{J}_n)\hat{g}''_n(\hat{u}_n)(\hat{u}'_n,\hat{u}'_n)\\
&+h^4\varphi_4(h\hat{J}_n)M(\hat{u}_n)+\mathcal{O}(h^5),
\end{split}
\end{equation}
where we denote by
\begin{equation}\label{Mun}
M(\hat{u}_n)=3\hat{g}_n''(\hat{u}_n)(\hat{u}'_n,\hat{u}''_n)+\hat{g}'''_n(\hat{u}_n)(\hat{u}'_n,\hat{u}'_n,\hat{u}'_n).
\end{equation}

\subsection{Stiff order conditions}
Here, we give the order conditions of methods  (\ref{generalDPG2}) and (\ref{generalDPG3}) to obtain a local truncation error of orders three and four, respectively. We define the local truncation error as 
$$\hat{e}_{n+1}=\bar{u}_{n+1}-\hat{u}_{n+1},$$
where $\hat{u}_{n+1}=u(t_{n+1})$ as in the previous section and $\bar{u}_{n+1}$ is the numerical solution obtained in (\ref{ExpEuler})-(\ref{generalDPG3}) considering $\hat{u}_{n}$ as an initial condition. Similarly, we denote by $\bar{u}_{n_2}$ and $\bar{u}_{n_3}$ the internal stages computed in (\ref{generalDPG2}) and (\ref{generalDPG3}) with $\hat{u}_{n}$ on the right-hand side.

\begin{proposition}
The hybrid Euler method defined in (\ref{ExpEuler}) is second order, i.e., 
$$\hat{e}_{n+1}=\mathcal{O}(h^3).$$
\begin{proof}
From (\ref{ExpEuler}) we have that
$$\bar{u}_{n+1}=\hat{u}_{n}+h^2\hat{J}_n\varphi_2(h\hat{J}_n)(\hat{g}_n(\hat{u}_{n})+\hat{J}_n\hat{u}_{n})+h(\hat{g}_n(\hat{u}_{n})+\hat{J}_n\hat{u}_{n}),$$
and from the recurrence formula (\ref{Recurrence})
\begin{equation}\label{ExpEuler2}
\bar{u}_{n+1}=\hat{u}_{n}+h\varphi_1(h\hat{J_n})(\hat{g}_n(\hat{u}_n)+\hat{J}_n\hat{u}_n),
\end{equation}
which is the classical exponential Euler method. Finally, subtracting (\ref{ExpEuler2}) from (\ref{TaylorExact}), we obtain the stated result. 
\end{proof}
\end{proposition}

\begin{theorem}\label{LocalTrunErrDPG2}
The two-stage DPG method defined in (\ref{generalDPG2}) with the coefficient functions satisfying
\begin{equation}\label{OrdCondDPG2}
\displaystyle{\left\{
\begin{split}
b_1(hJ_n)+b_2(hJ_n)&=\varphi_1(hJ_n),\\
\frac{1}{8}b_2(hJ_n)&=\varphi_3(hJ_n).\\
\end{split}
\right.} 
\end{equation}
is third order, i.e., 
$$\hat{e}_{n+1}=\mathcal{O}(h^4).$$
\begin{proof}
First, we need to verify that method (\ref{generalDPG2}) is stable for equilibrium points, i.e., it needs to reproduce constant solutions \footnote{From (\ref{TaylorExact}), it is easy to see that the numerical scheme should reproduce constant solutions to be at least first-order accurate. In classical exponential Runge-Kutta methods, enforcing the scheme to reproduce constants leads to order conditions for the coefficient functions in the internal stages. In the DPG methods we present here, this is satisfied by construction. We include this step here as a sanity check.}. If we suppose that $u\equiv \hat{u}_n$ then $\hat{g}_n(\hat{u}_n)=-\hat{J}_n\hat{u}_n$. We easily see that the second equation in (\ref{generalDPG2}) reproduces constant solutions by construction as in this case  $\bar{u}_{n_{2}}=\hat{u}_n$. We need to now enforce $\bar{u}_{n+1}=\hat{u}_n$, so from the first equation in (\ref{generalDPG2}), we obtain 
\begin{equation}\nonumber
\begin{split}
\hat{u}_n&=e^{h\hat{J}_n}\hat{u}_n-h\hat{J}_n\hat{u}_n(b_1(h\hat{J}_n)+b_2(h\hat{J}_n))\Longrightarrow(h\hat{J}_n)^{-1}(e^{h\hat{J}_n}-I)=b_1(h\hat{J}_n)+b_2(h\hat{J}_n)\\
&\Longrightarrow b_1(h\hat{J}_n)+b_2(h\hat{J}_n)=\varphi_1(h\hat{J}_n),
\end{split}
\end{equation}
which is the first condition in (\ref{OrdCondDPG2}). We now rewrite method (\ref{generalDPG2}) employing (\ref{Recurrence}) as
\begin{equation}\label{RewDPG2}
\displaystyle{\left\{
\begin{split}
\bar{u}_{n+1}&=\hat{u}_n+h\varphi_1(h\hat{J}_n)(\hat{g}_n(\hat{u}_{n})+\hat{J}_n\hat{u}_n)+hb_2(h\hat{J}_n)(\hat{g}_n(\bar{u}_{n_2})-\hat{g}_n(\hat{u}_n)),\\
\bar{u}_{n_2}&=\hat{u}_n+h\varphi_2(h\hat{J}_n)(\hat{g}_n(\hat{u}_{n})+\hat{J}_n\hat{u}_n).\\
\end{split}
\right.} 
\end{equation}
Next, we write the Taylor expansion of the numerical method to compare it with the analytical one (\ref{TaylorExact}). 
Employing the Taylor expansion of $\hat{g}_n$ at $\hat{u}_n$, we have that
$$\hat{g}_n(\bar{u}_{n_2})-\hat{g}_n(\hat{u}_n)=\frac{1}{2!}\hat{g}''_n(\hat{u}_n)(\bar{u}_{n_2}-\hat{u}_n,\bar{u}_{n_2}-\hat{u}_n)+\mathcal{R},$$
therefore, in the first equation of (\ref{RewDPG2}) we obtain 
\begin{equation}\label{TaylorNumDPG2}
\begin{split}
\bar{u}_{n+1}&=\hat{u}_n+h\varphi_1(hJ_n)(\hat{g}_n(\hat{u}_{n})+\hat{J}_n\hat{u}_n)+hb_2(h\hat{J}_n)\frac{1}{2!}\hat{g}''_n(\hat{u}_n)(\bar{u}_{n_2}-\hat{u}_n,\bar{u}_{n_2}-\hat{u}_n)+\mathcal{R}.\\
\end{split}
\end{equation}
We need to express $\bar{u}_{n_2}-\hat{u}_n$ in terms of $\hat{u}'_n=\hat{g}_n(\hat{u}_{n})+\hat{J}_n\hat{u}_n$. For that, we employ that
\begin{equation}\label{ChainRule}
\hat{J}_nu^{(k)}=u^{(k+1)}(t)-\frac{d^k}{dt^k}\hat{g}_n(u(t)),\;\forall k=0,1,\ldots;
\end{equation}
therefore, from the chain rule we have $\hat{J}_n\hat{u}'_n=\hat{u}''_n-\hat{g}'_n(\hat{u}_n)\hat{u}'_n=\hat{u}''_n$. From the second equation in (\ref{RewDPG2}) and (\ref{Recurrence}), we have that  
\begin{equation}\label{Un2DPG2}
\begin{split}
\bar{u}_{n_2}-\hat{u}_n&=h\varphi_2(h\hat{J}_n)\hat{u}'_n=h\left(\frac{1}{2!}+h\hat{J}_n\varphi_3(h\hat{J}_n)\right)\hat{u}'_n=\frac{h}{2!}\hat{u}'_n+h^2\varphi_3(h\hat{J}_n)\hat{u}''_n\\
&=\frac{h}{2!}\hat{u}'_n+h^2(\frac{1}{3!}+h\hat{J}_n\varphi_4(h\hat{J}_n))\hat{u}''_n=\frac{h}{2!}u'_n+\frac{h^2}{3!}u''_n+\mathcal{O}(h^3),\\
\end{split}
\end{equation}
where we have employed that matrices $\hat{J}_n$ and $\varphi_i(h\hat{J}_n)$ commute for all $i$. Finally, substituting (\ref{Un2DPG2}) into (\ref{RewDPG2})
$$\bar{u}_{n+1}=\hat{u}_n+h\varphi_1(hJ_n)(\hat{g}_n(\hat{u}_{n})+\hat{J}_n\hat{u}_n)+h^3b_2(h\hat{J}_n)\frac{1}{2!2!2!}\hat{g}''_n(\hat{u}_n)(u'_n,u'_n)+\mathcal{O}(h^4),$$
and subtracting from (\ref{TaylorExact}), we obtain that the local truncation error is of order $h^4$ provided that $b_2(h\hat{J_n})=8\varphi_3(h\hat{J}_n)$.
\end{proof}
\end{theorem}

\begin{theorem}
The three-stage DPG method defined in (\ref{generalDPG3}) with the coefficient functions  satisfying 
\begin{equation}\label{OrdCondDPG3}
\displaystyle{\left\{
\begin{split}
&b_1(hJ_n)+b_2(hJ_n)+b_3(hJ_n)=\varphi_1(hJ_n),\\
&\frac{1}{4}b_2(hJ_n)+b_3(hJ_n)=2\varphi_3(hJ_n),\\
&\frac{1}{8}b_2(hJ_n)+b_3(hJ_n)=6\varphi_4(hJ_n).\\
\end{split}
\right.} 
\end{equation}
and the correction factor $C_{n_2}:=-\frac{1}{4}g'_n(u_{n_2})(u_{n_3}-2u_{n_2}+u_n)$, is fourth order, i.e., 
$$\hat{e}_{n+1}=\mathcal{O}(h^5).$$
\begin{proof}
If we suppose that the solution is constant, i.e., $u\equiv \hat{u}_n$ and $\hat{g}_n(\hat{u}_n)=-\hat{J}_n\hat{u}_n$, then it is easy to verify that $\bar{u}_{n_3}=\bar{u}_{n_2}=\hat{u}_n$ in (\ref{generalDPG3}) by construction. Moreover, in this particular case, $\bar{C}_{n_2}=-\frac{1}{4}\hat{g}'_n(\bar{u}_{n_2})(\bar{u}_{n_3}-2\bar{u}_{n_2}+\hat{u}_n)=0$ and, in order to reproduce constant solutions in (\ref{generalDPG3}), we need to enforce that $b_1(h\hat{J}_n)+b_2(h\hat{J}_n)+b_3(h\hat{J}_n)=\varphi_1(h\hat{J}_n)$, which is the first condition in (\ref{OrdCondDPG3}). We now rewrite method (\ref{generalDPG2}) employing (\ref{Recurrence}) as
\begin{equation}\label{RewDPG3}
\displaystyle{\left\{
\begin{split}
\bar{u}_{n+1}&=\hat{u}_n+h\varphi_1(h\hat{J}_n)(\hat{g}_n(\hat{u}_{n})+\hat{J}_n\hat{u}_n)+hb_2(h\hat{J}_n)(\hat{g}_n(\bar{u}_{n_2})-\hat{g}_n(\hat{u}_n)+\bar{C}_{n_2})\\
&+hb_3(h\hat{J}_n)(\hat{g}_n(\bar{u}_{n_3})-\hat{g}_n(\hat{u}_n)),\\
\bar{u}_{n_2}&=\hat{u}_n+h\varphi_2(h\hat{J}_n)(\hat{g}_n(\hat{u}_{n})+\hat{J}_n\hat{u}_n),\\
\bar{u}_{n_3}&=\hat{u}_n+h\hat{J}_n\bar{u}_{n_2}+h\hat{g}_n(\hat{u}_n).\\
\end{split}
\right.} 
\end{equation}
We again employ the Taylor expansion of $\hat{g}_n$ and $\hat{u}_n$ as follows
$$\hat{g}_n(\bar{u}_{n_i})-\hat{g}_n(\hat{u}_n)=\frac{1}{2!}\hat{g}''_n(\hat{u}_n)(\bar{u}_{n_i}-\hat{u}_n,\bar{u}_{n_i}-\hat{u}_n)+\frac{1}{3!}\hat{g}'''_n(\hat{u}_n)(\bar{u}_{n_i}-\hat{u}_n,\bar{u}_{n_i}-\hat{u}_n,\bar{u}_{n_i}-\hat{u}_n)+\mathcal{R},$$
for $i=2,3,$ so the first equation in (\ref{RewDPG3}) now reads
\begin{equation}\label{TaylorNumDPG3}
\begin{split}
\bar{u}_{n+1}&=\hat{u}_n+h\varphi_1(h\hat{J}_n)(\hat{g}_n(\hat{u}_{n})+\hat{J}_n\hat{u}_n)\\
&+hb_2(h\hat{J}_n)\left(\frac{1}{2!}\hat{g}''_n(\hat{u}_n)(\bar{u}_{n_2}-\hat{u}_n,\bar{u}_{n_2}-\hat{u}_n)\right.\\
&\left.+\frac{1}{3!}\hat{g}'''_n(\hat{u}_n)(\bar{u}_{n_2}-\hat{u}_n,\bar{u}_{n_2}-\hat{u}_n,\bar{u}_{n_2}-\hat{u}_n)+\bar{C}_{n_2}\right)\\
&+hb_3(h\hat{J}_n)\left(\frac{1}{2!}\hat{g}''_n(\hat{u}_n)(\bar{u}_{n_3}-\hat{u}_n,\bar{u}_{n_3}-\hat{u}_n)\right.\\
&\left.+\frac{1}{3!}\hat{g}'''_n(\hat{u}_n)(\bar{u}_{n_3}-\hat{u}_n,\bar{u}_{n_3}-\hat{u}_n,\bar{u}_{n_3}-\hat{u}_n)\right)+\mathcal{R}.
\end{split}
\end{equation}

We employ the following identities derived from (\ref{ChainRule}) by applying the chain rule
\begin{equation}\nonumber
\displaystyle{\left\{
\begin{split}
\hat{J}_n\hat{u}'_n&=\hat{u}''_n-\hat{g}'_n(\hat{u}_n)\hat{u}'_n=\hat{u}''_n,\\
\hat{J}_n\hat{u}''_n&=\hat{u}'''_n-\hat{g}''_n(\hat{u}_n)(\hat{u}'_n,\hat{u}'_n)-\hat{g}'_n(\hat{u}_n)\hat{u}''_n=\hat{u}'''_n-\hat{g}''_n(\hat{u}_n)(\hat{u}'_n,\hat{u}'_n),\\
\end{split}
\right.} 
\end{equation}
and also that $\hat{u}'_n=\hat{g}_n(\hat{u}_{n})+\hat{J}_n\hat{u}_n$. From the second equation in (\ref{RewDPG3}) and (\ref{Recurrence}), we obtain
\begin{equation}\label{Un2DPG3}
\begin{split}
\bar{u}_{n_2}-\hat{u}_n&=h\varphi_2(hJ_n)u'_n=h\left(\frac{1}{2!}+h\hat{J}_n\varphi_3(h\hat{J}_n)\right)\hat{u}'_n=\frac{h}{2!}\hat{u}'_n+h^2\varphi_3(h\hat{J}_n)\hat{u}''_n\\
&=\frac{h}{2!}\hat{u}'_n+h^2\left(\frac{1}{3!}+h\hat{J}_n\varphi_4(h\hat{J}_n)\right)\hat{u}''_n\\
&=\frac{h}{2!}\hat{u}'_n+\frac{h^2}{3!}\hat{u}''_n+h^3\varphi_4(h\hat{J}_n)\left(\hat{u}'''_n-\hat{g}''_n(\hat{u}_n)(\hat{u}'_n,\hat{u}'_n)\right),\\
\end{split}
\end{equation}
and from the third equation of (\ref{RewDPG3})
\begin{equation}\label{Un3DPG3}
\begin{split}
\bar{u}_{n_3}-\hat{u}_n&=h\hat{J}_n\bar{u}_{n_2}+h\hat{g}_n(\hat{u}_n)=h\hat{J}_n\left(\hat{u}_n+h\varphi_2(h\hat{J}_n)\hat{u}'_n\right)+h\hat{g}_n(\hat{u}_n)\\
&=h\hat{u}'_n+h^2\hat{J}_n\varphi_2(h\hat{J}_n)\hat{u}'_n=h\hat{u}'_n+h^2\left(\frac{1}{2!}+h\hat{J}_n\varphi_3(h\hat{J}_n)\right)\hat{u}''_n\\
&=h\hat{u}'_n+\frac{h^2}{2!}\hat{u}''_n+h^3\varphi_3(h\hat{J}_n)(\hat{u}'''_n-\hat{g}''_n(\hat{u}_n)(\hat{u}'_n,\hat{u}'_n)).\\
\end{split}
\end{equation}
Finally, substituting (\ref{Un2DPG3}) and (\ref{Un3DPG3}) in (\ref{TaylorNumDPG3}), we obtain 
\begin{equation}\label{TaylorNumDPG3New}
\begin{split}
\bar{u}_{n+1}&=\hat{u}_n+h\varphi_1(h\hat{J}_n)(\hat{g}_n(\hat{u}_{n})+\hat{J}_n\hat{u}_n)\\
&+hb_2(h\hat{J}_n)\left(\frac{h^2}{2!2!2!}\hat{g}''_n(\hat{u}_n)(\hat{u}'_n,\hat{u}'_n)+\frac{h^3}{2!3!}\hat{g}''_n(\hat{u}_n)(\hat{u}''_n,\hat{u}'_n)\right.\\
&+\left.\frac{h^3}{3!2!2!2!}\hat{g}'''_n(\hat{u}_n)(\hat{u}'_n,\hat{u}'_n,\hat{u}'_n)+\bar{C}_{n_2}\right)\\
&+hb_3(h\hat{J}_n)\left(\frac{h^2}{2!}\hat{g}''_n(\hat{u}_n)(\hat{u}'_n,\hat{u}'_n)+\frac{h^3}{2!}\hat{g}''_n(\hat{u}_n)(\hat{u}''_n,\hat{u}'_n)\right.\\
&+\left.\frac{h^3}{3!}\hat{g}'''_n(\hat{u}_n)(\hat{u}'_n,\hat{u}'_n,\hat{u}'_n)\right)+\mathcal{O}(h^5).\\
\end{split}
\end{equation}

We employ $\bar{C}_{n_2}$ to correct the coefficients in (\ref{TaylorNumDPG3New}) and factor out the term $M(\hat{u}_n)$ defined in (\ref{Mun}). First, note that from (\ref{Un2DPG3}) and (\ref{Un3DPG3}), we have that 
\begin{equation}\label{OrderVarespilon}
\bar{u}_{n_3}-2\bar{u}_{n_2}+\hat{u}_n=\frac{h^2}{3!}\hat{u}''_{n}+\mathcal{O}(h^3),\;\;\;\bar{u}_{n_2}-\hat{u}_n=\frac{h}{2!}\hat{u}'_n+\frac{h^2}{3!}\hat{u}''_n+\mathcal{O}(h^3),
\end{equation}
and employing Taylor expansions 
\begin{equation}
\begin{split}
\bar{C}_{n_2}=&-\frac{1}{4}\hat{g}'_n(\bar{u}_{n_2})(\bar{u}_{n_3}-2\bar{u}_{n_2}+\hat{u}_n)=-\frac{1}{4}\frac{h^2}{3!}\hat{g}'_n(\bar{u}_{n_2})\hat{u}''_n+\mathcal{O}(h^3)\\
=&-\frac{1}{4}\frac{h^2}{3!}\hat{g}'_n(\hat{u}_n)\hat{u}''_n-\frac{1}{4}\frac{h^2}{3!}\hat{g}''(\hat{u}_n)(\bar{u}_{n_2}-\hat{u}_n,\hat{u}''_n)\\
&-\frac{1}{4}\frac{h^2}{3!}\hat{g}'''(\hat{u}_n)(\bar{u}_{n_2}-\hat{u}_n,\bar{u}_{n_2}-\hat{u}_n,\hat{u}''_n)+\mathcal{O}(h^3)\\
=&-\frac{1}{4}\frac{h^3}{2!3!}g''(\hat{u}_n)(\hat{u}'_n,\hat{u}''_n)+\mathcal{O}(h^4). 
\end{split}
\end{equation}
Therefore, we can express (\ref{TaylorNumDPG3New}) in terms of $M(\hat{u}_n)$ and subtracting it from (\ref{TaylorExact}), we obtain that the local truncation error is of order $h^5$ provided that the coefficients $b_2(h\hat{J}_n)$ and $b_3(h\hat{J}_n)$ satisfy (\ref{OrdCondDPG3}). 
\end{proof}
\end{theorem}

\subsection{Final methods.}
In Table \ref{FinalCoeff}, we summarize the coefficients for the updates in methods (\ref{generalDPG2}) and (\ref{generalDPG3}). 

\renewcommand{\arraystretch}{1.8}
\begin{table}[h!]
\begin{adjustbox}{width=1\textwidth}
\centering
\begin{tabular}{|c|c|c|c|}\hline
method&$b_1(hJ_n)$&$b_2(hJ_n)$&$b_3(hJ_n)$\\ \hline
two-stage DPG&$\varphi_1(hJ_n)-8\varphi_3(hJ_n)$&$8\varphi_3(hJ_n)$&$-$\\ \hline
three-stage DPG&$\varphi_1(hJ_n)-14\varphi_3(hJ_n)+36\varphi_4(hJ_n)$&$16\varphi_3(hJ_n)-48\varphi_4(hJ_n)$&$12\varphi_4(hJ_n)-2\varphi_3(hJ_n)$\\ \hline
\end{tabular}
\end{adjustbox}
\caption{Coefficients of the two-stage (\ref{generalDPG2}) and three-stage (\ref{generalDPG3}) DPG methods.}
\label{FinalCoeff}
\end{table}

Note that the second and third-order methods (\ref{ExpEuler}) and (\ref{generalDPG2}) have the same cost as classical exponential Rosenbrock methods of one and two stages, respectively. However, the fourth-order method (\ref{generalDPG3}) requires one less exponential step than the classical three-stage methods, being therefore cheaper. Finally, the extension to variable time-step size of the methods presented here is straightforward; we have considered a uniform partition for simplicity. Moreover, methods (\ref{generalDPG2}) and (\ref{generalDPG3}) are nested because the internal stage that requires an exponential step is the same in both methods so they can be employed to adapt the time step size.

\section{Convergence}\label{Sec:Convergence}
In this section, we prove the convergence of the methods described in Section \ref{Sec:Methods}. We focus on methods (\ref{generalDPG2}) and (\ref{generalDPG3}) as method (\ref{ExpEuler}) is equivalent to the standard exponential Euler method, and its convergence proof is stated in \cite{hochbruck2009exponential,luan2014exponential}. We will follow the strategy presented in \cite{luan2014exponential}, and for that, we define the following errors, 
\begin{equation}
\begin{split}
e_{n+1}=u_{n+1}-\hat{u}_{n+1},&\;\;\bar{e}_{n+1}=u_{n+1}-\bar{u}_{n+1},\\
\bar{E}_{n_2}=u_{n_2}-\bar{u}_{n_2},&\;\;\bar{E}_{n_3}=u_{n_3}-\bar{u}_{n_3}.
\end{split}
\end{equation}
We need to bound the total error $e_{n+1}$. Note that $e_{n+1}=\bar{e}_{n+1}+\hat{e}_{n+1}$, where $\hat{e}_{n+1}$ is the local truncation error we estimated in Section \ref{Sec:Stiff} (See \ref{App:Summarize} for a summary of notation.)

\subsection{Preliminary error bounds and stability assumptions}
In order to prove the final convergence, we need the preliminary error bounds that are presented in the following lemma. 
\begin{lemma}\label{LemmaPre}
Under Assumptions \ref{Ass1} and \ref{Ass2}, there exists a constant $C\geq0$ such that the following estimates hold
\begin{subequations}
\begin{equation}\label{pre1}
\left\lVert g'_n(\hat{u}_n)\right\rVert_{_{\mathcal{L}(X)}}=\Vert J(\hat{u}_n)-J_n\Vert_{_{\mathcal{L}(X)}} \leq C\Vert e_n\Vert ,
\end{equation}
\begin{equation}\label{pre2}
\Vert g_n(u_n)-g_n(\hat{u}_n)\Vert \leq C\Vert e_n\Vert ^2,
\end{equation}
\begin{equation}\label{pre3}
\Vert \varphi_j(hJ_n)-\varphi_j(h\hat{J}_n)\Vert_{_{\mathcal{L}(X)}} \leq Ch\Vert e_n\Vert ,
\end{equation}
\begin{equation}\label{pre4}
\Vert b_j(hJ_n)-b_j(h\hat{J}_n)\Vert_{_{\mathcal{L}(X)}} \leq Ch\Vert e_n\Vert ,
\end{equation}
\end{subequations}
for all $i$ and $j$. 
\begin{proof}
The first estimate (\ref{pre1}) is a direct consequence of the linearization (\ref{JacRem}) and the Lipschitz condition given in Assumption \ref{Ass1}. Estimate (\ref{pre2}) is obtained from (\ref{pre1}) and the following Taylor expansion
$$g_n(u_n)-g_n(\hat{u}_n)=g'_n(\hat{u}_n)e_n+\mathcal{O}(\Vert e_n\Vert ^2).$$
Finally, (\ref{pre3}) and (\ref{pre4}) are obtained employing the definition of the $\varphi-$functions, the fact that the coefficients $b_j(hJ_n)$ are linear combinations of $\varphi-$functions, and inequalities (\ref{Ass1for}) and (\ref{Ass2for}) (see Lemmas 4.2 and 4.3 in \cite{luan2014exponential}).
\end{proof}
\end{lemma}

Finally, as $J_n$ varies from time step to time step, we will also employ the following stability assumption from \cite{luan2014exponential}. 
\begin{assumption}\label{Ass3}
We assume that the following stability condition holds
\begin{equation}
\left\lVert\prod_{k=0}^{n-\nu} e^{hJ_{n-k}} \right\rVert_{\mathcal{L}(X)}\leq C_S,\;\;\; 0\leq t_\nu \leq t_n\leq T,
\end{equation}
where $C_S$ is a constant that is uniform in $\nu$ and $n$. 
\end{assumption}

\begin{remark}
Both in the construction of the methods and the convergence proofs we present in this section, we consider a uniform time step for simplicity. However, the extension to a variable time-step size is straightforward under the usual mild restrictions on the step size sequence (see Section 4 in \cite{hochbruck2009exponential} for details). 
\end{remark}

\subsection{Convergence of the two-stage DPG method}
Subtracting (\ref{RewDPG2}) from (\ref{generalDPG2}) and employing the identity
$$F(u)=g_n(u)+J_n(u)=\hat{g}_n(u)+\hat{J}_n(u),$$
we have that 
\begin{equation}\label{errorDPG2}
\begin{split}
\bar{e}_{n+1}&=e^{hJ_n}e_n+h\varphi_1(hJ_n)(g_n(u_{n})-g_n(\hat{u}_n))+h(\varphi_1(hJ_n)-\varphi_1(h\hat{J}_n))F(\hat{u}_n)\\
&+hb_2(hJ_n)(D_{n_2}-\hat{D}_{n_2})+h(b_2(hJ_n)-b_2(h\hat{J}_n))\hat{D}_{n_2}, 
\end{split}
\end{equation}
where we denote by
$$D_{n_2}:=g_n(u_{n_2})-g_n(u_n),\;\;\hat{D}_{n_2}:=\hat{g}_n(\bar{u}_{n_2})-\hat{g}_n(\hat{u}_n).$$

Therefore, form (\ref{errorDPG2}) we obtain 
\begin{equation}\label{Recursion}
e_{n+1}=e^{hJ_n}e_n+hq_n+hQ_{n}+\hat{e}_{n+1},
\end{equation}
where 
$$q_n=\varphi_1(hJ_n)(g_n(u_{n})-g_n(\hat{u}_n))+(\varphi_1(hJ_n)-\varphi_1(h\hat{J}_n))F(\hat{u}_n),$$
$$Q_{n}=b_2(hJ_n)(D_{n_2}-\hat{D}_{n_2})+(b_2(hJ_n)-b_2(h\hat{J}_n))\hat{D}_{n_2}.$$

\begin{proposition}\label{ErrsDPG2}
Under Assumptions (\ref{Ass1}) and (\ref{Ass2}) we have that 
\begin{subequations}
\begin{equation}\label{Err1DPG2}
\Vert q_n\Vert \leq Ch\Vert e_n\Vert +C\Vert e_n\Vert ^2,
\end{equation}
\begin{equation}\label{Err2DPG2}
\Vert Q_{n}\Vert \leq  Ch\Vert e_n\Vert +C\Vert e_n\Vert ^2+C(h+\Vert e_n\Vert +\Vert \bar{E}_{n_2}\Vert )\Vert \bar{E}_{n_2}\Vert ,
\end{equation}
\begin{equation}\label{Err3DPG2}
\Vert \bar{E}_{n_2}\Vert \leq  C\Vert e_n\Vert +Ch||e_n||^2.
\end{equation}
\end{subequations} 
\begin{proof}
Estimate (\ref{Err1DPG2}) is a direct consequence of (\ref{pre2}) and (\ref{pre3}). To obtain (\ref{Err1DPG2}) we need to bound $D_{n_2}-\hat{D}_{n_2}$. For that, we rewrite
$$D_{n_2}-\hat{D}_{n_2}=(g_n(u_{n_2})-g_n(\bar{u}_{n_2}))+(g_n(\bar{u}_{n_2})-\hat{g}_n(\bar{u}_{n_2}))+(\hat{g}_n(\hat{u}_{n})-g_n(\hat{u}_{n}))+(g_n(\hat{u}_{n})-g_n(u_{n})),$$
and we bound each term separately. We employ the Taylor expansion
$$g_n(u_{n_2})-g_n(\bar{u}_{n_2})=g'_n(\bar{u}_{n_2})\bar{E}_{n_2}+\mathcal{O}(\Vert \bar{E}_{n_2}\Vert ^2)=(J(\bar{u}_{n_2})-J_n)\bar{E}_{n_2}+\mathcal{O}(\Vert \bar{E}_{n_2}\Vert ^2),$$
and also the identity 
$$g_n(u)-\hat{g}_n(u)=F(u)-J_nu-F(u)+\hat{J}_nu=(\hat{J}_n-J_n)u,$$
so we have 
$$g_n(\bar{u}_{n_2})-\hat{g}_n(\bar{u}_{n_2})=(\hat{J}_n-J_n)\bar{u}_{n_2},\;\;\hat{g}_n(\hat{u}_{n})-g_n(\hat{u}_{n})=(J_n-\hat{J}_n)\hat{u}_{n},$$
and therefore, 
\begin{equation}\nonumber
D_{n_2}-\hat{D}_{n_2}=(J(\bar{u}_{n_2})-J_n)\bar{E}_{n_2}+(\hat{J}_n-J_n)\bar{u}_{n_2}+(J_n-\hat{J}_n)\hat{u}_{n}+g_n(\hat{u}_n)-g_n(u_n)+\mathcal{O}(\Vert \bar{E}_{n_2}\Vert ^2).
\end{equation}
From the Lipschitz condition (\ref{Ass1for}) and (\ref{pre3}) we get 
$$\Vert D_{n_2}-\hat{D}_{n_2}\Vert\leq C\Vert\bar{u}_{n_2}-u_n\Vert\Vert\bar{E}_{n_2}\Vert+C\Vert e_n\Vert\Vert\bar{u}_{n_2}-\hat{u}_n\Vert+C\Vert e_n\Vert^2+C\Vert\bar{E}_{n_2}\Vert^2,$$
from (\ref{Un2DPG2}) we know that $\bar{u}_{n_2}-\hat{u}_n=\mathcal{O}(h)$, therefore
$$\Vert \bar{u}_{n_2}-\hat{u}_n+\hat{u}_n-u_n\Vert\leq Ch+C\Vert e_n\Vert,$$
and employing (\ref{pre4}), we obtain estimate (\ref{Err2DPG2}).
For the last estimate, we have that
$$\bar{E}_{n_2}=e_n+h\varphi_2(hJ_n)(F(u_n)-F(\hat{u}_n))+h(\varphi_2(hJ_n)-\varphi_2(h\hat{J}_n))F(\hat{u}_n).$$
Employing the identity
$$F(u_n)-F(\hat{u}_n)=J_ne_n+g_n(u_n)-g_n(\hat{u}_n),$$
and recurrence formula (\ref{Recurrence}), we obtain
$$\bar{E}_{n_2}=\varphi_1(hJ_n)e_n+h\varphi_2(hJ_n)(g_n(u_n)-g_n(\hat{u}_n))+h(\varphi_2(hJ_n)-\varphi_2(h\hat{J}_n))F(\hat{u}_n).$$
Finally, from Lemma \ref{LemmaPre} we obtain (\ref{Err3DPG2}).
\end{proof}
\end{proposition}

We provide the final convergence result for the two-stage DPG method in the following theorem.

\begin{theorem}\label{TheoDPG2}
Under Assumptions \ref{Ass1}, \ref{Ass2}, and \ref{Ass3}, we conclude that the two-stage DPG method (\ref{generalDPG2}) with coefficients satisfying (\ref{OrdCondDPG2}), is third-order, i.e., the numerical solution satisfies the error bound
\begin{equation}\label{FinalErrorDPG2}
||e_n||\leq Ch^{4}
\end{equation}
uniformly in $t_n$ and with constant $C$ independent of the time-step size. 
\begin{proof}
Solving recursion (\ref{Recursion}) with $e_0=0$, we obtain
$$e_n=\sum_{j=0}^{n-1}h\prod_{k=1}^{n-j-1}e^{hJ_{n-k}}(q_j+Q_j+h^{-1}\hat{e}_{j+1}).$$
Employing Proposition (\ref{ErrsDPG2}) and Theorem (\ref{LocalTrunErrDPG2}) we have the bound
$$||q_j||+||Q_j||+h^{-1}||\hat{e}_{j+1}||\leq C(h||e_j||+||e_j||^2+h^{3}),$$
and employing Assumption \ref{Ass3} we have 
$$||e_n||\leq C\sum_{j=0}^{n-1}h(h||e_j||+||e_j||^2+h^{3}).$$
Finally, the application of a discrete Gronwall lemma (see \cite{hochbruck2010exponential}, Lemma 2.15) yields the final error bound (\ref{FinalErrorDPG2}).
\end{proof}
\end{theorem}

\subsection{Convergence of the three-stage DPG method}
Subtracting (\ref{RewDPG3}) from (\ref{generalDPG3}), we have that 
\begin{equation}\label{errorDPG3}
\begin{split}
\bar{e}_{n+1}&=e^{hJ_n}e_n+h\varphi_1(hJ_n)(g_n(u_{n})-g_n(\hat{u}_n))+h(\varphi_1(hJ_n)-\varphi_1(h\hat{J}_n))F(\hat{u}_n)\\
&+hb_2(hJ_n)(D_{n_2}+C_{n_2}-\hat{D}_{n_2}-\bar{C}_{n_2})+h(b_2(hJ_n)-b_2(h\hat{J}_n))(\hat{D}_{n_2}+\bar{C}_{n_2}),\\ 
&+hb_3(hJ_n)(D_{n_3}-\hat{D}_{n_3})+h(b_3(hJ_n)-b_3(h\hat{J}_n))\hat{D}_{n_3}, 
\end{split}
\end{equation}
where
$$D_{n_i}=g_n(u_{n_i})-g_n(u_n),\;\;\hat{D}_{n_i}=\hat{g}_n(\bar{u}_{n_i})-\hat{g}_n(\hat{u}_n),\;\;i=2,3.$$

Therefore, form (\ref{errorDPG2}) we have that 
$$e_{n+1}=e^{hJ_n}e_n+hq_n+hQ_{n}+\hat{e}_{n+1},$$
where 
$$q_n=\varphi_1(hJ_n)(g_n(u_{n})-g_n(\hat{u}_n))+(\varphi_1(hJ_n)-\varphi_1(h\hat{J}_n))F(\hat{u}_n),$$
\begin{equation}\nonumber
\begin{split}
Q_{n}&=b_2(hJ_n)(D_{n_2}+C_{n_2}-\hat{D}_{n_2}-\bar{C}_{n_2})+(b_2(hJ_n)-b_2(h\hat{J}_n))(\hat{D}_{n_2}+\bar{C}_{n_2})\\
&+b_3(hJ_n)(D_{n_3}-\hat{D}_{n_3})+h(b_3(hJ_n)-b_3(h\hat{J}_n))\hat{D}_{n_3}.
\end{split}
\end{equation}

\begin{proposition}\label{ErrsDPG3}
Under Assumptions (\ref{Ass1}) and (\ref{Ass2}) we have that 
\begin{subequations}
\begin{equation}\label{Err1DPG3}
||q_n||\leq Ch||e_n||+C||e_n||^2,
\end{equation}
\begin{equation}\label{Err2DPG3}
||Q_{n}||\leq  Ch||e_n||+C||e_n||^2+C\sum_{j=2}^3(h+||e_n||+||\bar{E}_{n_j}||)||\bar{E}_{n_j}||,
\end{equation}
\begin{equation}\label{Err3DPG3}
\Vert \bar{E}_{n_j}\Vert \leq  C\Vert e_n\Vert +Ch||e_n||^2,\;\;j=2,3.
\end{equation}
\end{subequations} 
\begin{proof}
Estimate (\ref{Err1DPG3}) is the same as (\ref{Err1DPG2}). To obtain estimate (\ref{Err2DPG3}), we proceed as in (\ref{Err2DPG2}) to bound $D_{n_i}-\hat{D}_{n_i},\;i=2,3$, and additionally, we have to bound $C_{n_2}-\bar{C}_{n_2}$. For that, we rewrite
$$C_{n_2}-\bar{C}_{n_2}=\frac{1}{4}g'_n(u_{n_2})(\varepsilon_n-\bar{\varepsilon}_n)+\frac{1}{4}(\hat{g}'_n(\bar{u}_{n_2})-g'_n(u_{n_2}))\bar{\varepsilon}_n,$$
where we have denoted 
$$\varepsilon_n:=u_{n_3}-2u_{n_2}+u_n,\;\;\bar{\varepsilon}_n:=\bar{u}_{n_3}-2\bar{u}_{n_2}+\bar{u}_n.$$
We employ that $\varepsilon_n-\bar{\varepsilon}_n=E_{n_3}-2E_{n_2}+e_n$ and from (\ref{OrderVarespilon}), we know that $\bar{\varepsilon}_n=\mathcal{O}(h^2)$. Now, from the following estimates 
$$||g'_n(u_{n_2})||_{_{\mathcal{L}(X)}}=||J(u_{n_2})-J_n||_{_{\mathcal{L}(X)}}\leq C||u_{n_2}-u_n||\leq Ch,$$
$$||\hat{g}'_n(\bar{u}_{n_2})-g'_n(u_{n_2})||_{_{\mathcal{L}(X)}}=||J(u_{n_2})-J_n-J(\bar{u}_{n_2})+\hat{J}_n||_{_{\mathcal{L}(X)}}\leq C||E_{n_2}||+C||e_n||,$$
we obtain that 
$$||C_{n_2}-\bar{C}_{n_2}||\leq Ch||e_n||+Ch^2||e_n||+Ch||E_{n_2}||+Ch^2||E_{n_2}||+Ch||E_{n_3}||.$$
Combining this last estimate with the bounds for $D_{n_i}-\hat{D}_{n_i}$, we derive (\ref{Err2DPG3}). Finally, estimate (\ref{Err3DPG3}) for $j=2$ is the same as (\ref{Err3DPG2}). For $j=3$, we have
$$\bar{E}_{n_3}=e^{hJ_n}e_n+h\varphi_1(hJ_n)(g_n(u_n)-g_n(\hat{u}_n))+h(\varphi_1(hJ_n)-\varphi_1(h\hat{J}_n))F(\hat{u}_n),$$
and from Lemma \ref{LemmaPre} we obtain (\ref{Err3DPG3}).
\end{proof}
\end{proposition}

\begin{theorem}\label{TheoDPG3}
Under Assumptions \ref{Ass1}, \ref{Ass2}, and \ref{Ass3}, we conclude that the three-stage DPG method (\ref{generalDPG3}) with coefficients satisfying (\ref{OrdCondDPG3}) is fourth-order, i.e., the numerical solution satisfies the error bound
\begin{equation}\label{FinalErrorDPG2}
||e_n||\leq Ch^{5},
\end{equation}
uniformly in $t_n$ and with constant $C$ independent of the time-step size. 
\begin{proof}
Employing the preliminary error bounds from Proposition \ref{ErrsDPG3}, the final convergence proof is analogue to the one presented in Theorem \ref{TheoDPG2}.
\end{proof}
\end{theorem}

\section{Numerical results}\label{Sec:Results}
In this section, we first test the convergence predicted by the theory in a 2D+time semilinear equation. Additionally, we numerically observe that our methods satisfy the energy decaying property in a dissipative model and the Kinetic energy preservation in a conservative model.
\subsection{Hochbruck-Ostermann equation}
We consider the Hochbruck-Ostermann equation \cite{hochbruck2005explicit} in $\Omega=(0,1)^2$
\begin{equation}\label{HochOster}
u_t-\Delta u=\frac{1}{1+u^2}+f(x,y,t),
\end{equation}
with homogeneous Dirichlet boundary conditions and final time $T=1$. We select the linear source $f$ and the initial condition $u_0$ in such a way that the exact solution is $u(x,y,t)=x(1-x)y(1-y)e^t$. We discretize (\ref{HochOster}) in space by standard second-order central finite differences with a fixed grid of $2^6$ points in each dimension. 

\begin{figure}[h!]
\centering
\begin{tikzpicture}
\begin{axis}[small,
	width=1.75*\width,
	height=1.25*\height,
	ymode=log,
	xmode=log,
	xlabel=$h$,
	ylabel=Error at final time,
	grid=both,
    legend style={font=\small,at={(0.5,1.03)},anchor=south},
    legend entries={Hybrid Euler,DPG2,DPG3,Exponential Euler,expbr32,expbr43,Rate 2, Rate 3, Rate 4},
    legend columns=3
	]      
\addplot[
	color=black,dotted,line width=1pt, 
	mark=*,mark size=1.5pt,    
	]
table[x=tau,y=Err_EulerExpmv]{Results/ErrHochOster_Ros_r6_m10.txt};
\addplot[
	color=blue,dotted,line width=1pt, 
	mark=*,mark size=1.5pt,    
	]
table[x=tau,y=Err_RK2Expmv]{Results/ErrHochOster_Ros_r6_m10.txt};
\addplot[
	color=red,dotted,line width=1pt, 
	mark=*,mark size=1.5pt,    
	]
table[x=tau,y=Err_RK3Expmv]{Results/ErrHochOster_Ros_r6_m10.txt};
\addplot[
	color=black,solid,line width=0.6pt, 
	mark=*,mark size=1.5pt,    
	]
table[x=tau,y=Err_EulerExpmv]{Results/ErrHochOster_Class_Ros_r6_m10.txt};
\addplot[
	color=blue,solid,line width=0.6pt, 
	mark=*,mark size=1.5pt,    
	]
table[x=tau,y=Err_RK2Expmv]{Results/ErrHochOster_Class_Ros_r6_m10.txt};
\addplot[
	color=red,solid,line width=0.6pt, 
	mark=*,mark size=1.5pt,    
	]
table[x=tau,y=Err_RK3Expmv]{Results/ErrHochOster_Class_Ros_r6_m10.txt};
\addplot[
	color=black,solid,line width=0.6pt, 
	]
table[x=tau,y=Slope]{Results/Slope2.txt};
\addplot[
	color=blue,solid,line width=0.6pt, 
	]
table[x=tau,y=Slope]{Results/Slope3.txt};
\addplot[
	color=red,solid,line width=0.6pt, 
	]
table[x=tau,y=Slope]{Results/Slope4.txt};
\end{axis}
\end{tikzpicture}
\caption{Convergence in time at $T=1$ for the hybrid exponential Euler (\ref{ExpEuler}), the two-stage DPG (\ref{generalDPG2}), and the three-stage DPG (\ref{generalDPG3}) methods, and the classical exponential Euler, expbr32, and expbr32 methods.}
\label{ConvergenceHochOster}
\end{figure}
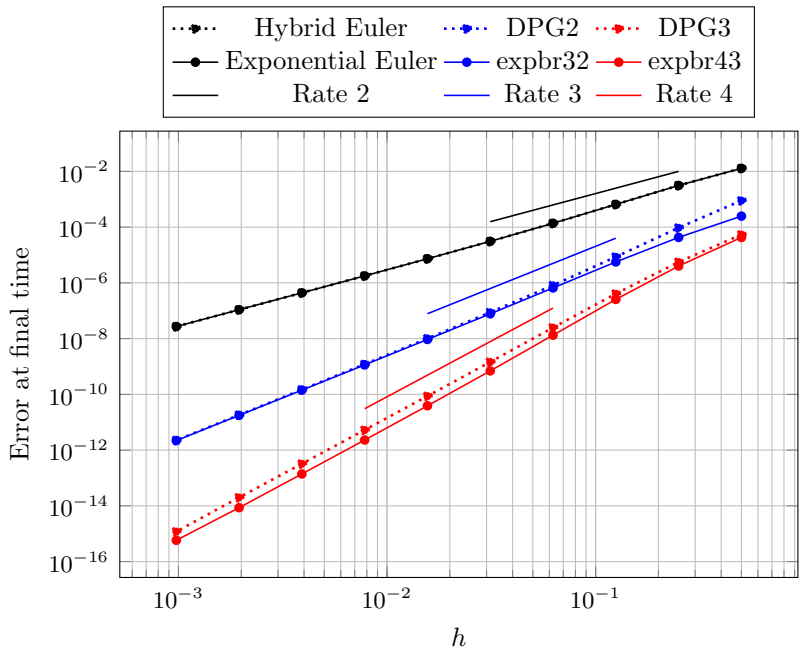

Figure \ref{ConvergenceHochOster} shows the convergence of the error in time by solving the resulting system of ODEs by the hybrid exponential Euler method (\ref{ExpEuler}), the two-stage DPG method (\ref{generalDPG2}), and the three-stage DPG method (\ref{generalDPG3}). Note that this problem is semilinear and non-autonomous, so we implemented the DPG methods taking into account Remarks \ref{rmk:1} and \ref{rmk:2}. We measure the infinity norm at the final time and observe that the convergence rates are the ones predicted in Section \ref{Sec:Convergence}. In Figure \ref{ConvergenceHochOster}, we also compare the convergence of our three methods with the classical exponential Rosenbrock methods (exponential Euler, expbr32, expbr43) form \cite{hochbruck2009exponential}. We observe that both Euler methods are equivalent (as stated in Section \ref{Sec:Methods}) and conclude that the two- and three-stage DPG multistage methods are as accurate as the classical ones.

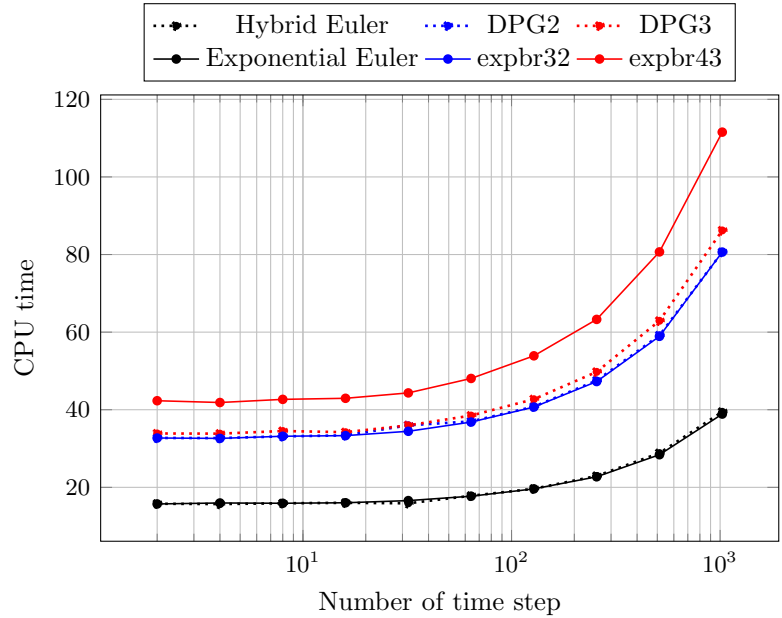
\begin{figure}[h!]
\centering
\begin{tikzpicture}
\begin{axis}[small,
	width=1.75*\width,
	height=1.25*\height,
	xmode=log,
	xlabel=Number of time step,
	ylabel=CPU time,
	grid=both,
    legend style={font=\small,at={(0.5,1.03)},anchor=south},
    legend entries={Hybrid Euler,DPG2,DPG3,Exponential Euler,expbr32,expbr43},
    legend columns=3
	]      
\addplot[
	color=black,dotted,line width=1pt, 
	mark=*,mark size=1.5pt,    
	]
table[x=Nsteps,y=Times_EulerExpmv]{Results/TimesHochOster_Ros_r100_m10.txt};
\addplot[
	color=blue,dotted,line width=1pt, 
	mark=*,mark size=1.5pt,    
	]
table[x=Nsteps,y=Times_RK2Expmv]{Results/TimesHochOster_Ros_r100_m10.txt};
\addplot[
	color=red,dotted,line width=1pt, 
	mark=*,mark size=1.5pt,    
	]
table[x=Nsteps,y=Times_RK3Expmv]{Results/TimesHochOster_Ros_r100_m10.txt};
\addplot[
	color=black,solid,line width=0.6pt, 
	mark=*,mark size=1.5pt,    
	]
table[x=Nsteps,y=Times_EulerExpmv]{Results/TimesHochOster_Class_Ros_r100_m10.txt};
\addplot[
	color=blue,solid,line width=0.6pt, 
	mark=*,mark size=1.5pt,    
	]
table[x=Nsteps,y=Times_RK2Expmv]{Results/TimesHochOster_Class_Ros_r100_m10.txt};
\addplot[
	color=red,solid,line width=0.6pt, 
	mark=*,mark size=1.5pt,    
	]
table[x=Nsteps,y=Times_RK3Expmv]{Results/TimesHochOster_Class_Ros_r100_m10.txt};
\end{axis}
\end{tikzpicture}
\caption{CPU time versus number of time steps for the hybrid exponential Euler (\ref{ExpEuler}), the two-stage DPG (\ref{generalDPG2}), and the three-stage DPG (\ref{generalDPG3}) methods, and the classical exponential Euler, expbr32, and expbr32 methods.}
\label{CostHochOster}
\end{figure}

Finally, in \ref{CostHochOster} we compare the cost of the multistage DPG methods with the classical exponential Rosenbrock methods form \cite{hochbruck2009exponential}. We conclude that that the fourth order DPG method has the same cost of the third order methods due to the postprocessing step in one of the internal stages.

\subsection{Allen-Cahn equation}

\begin{figure}[h!]
\center{
\includegraphics[scale=0.5]{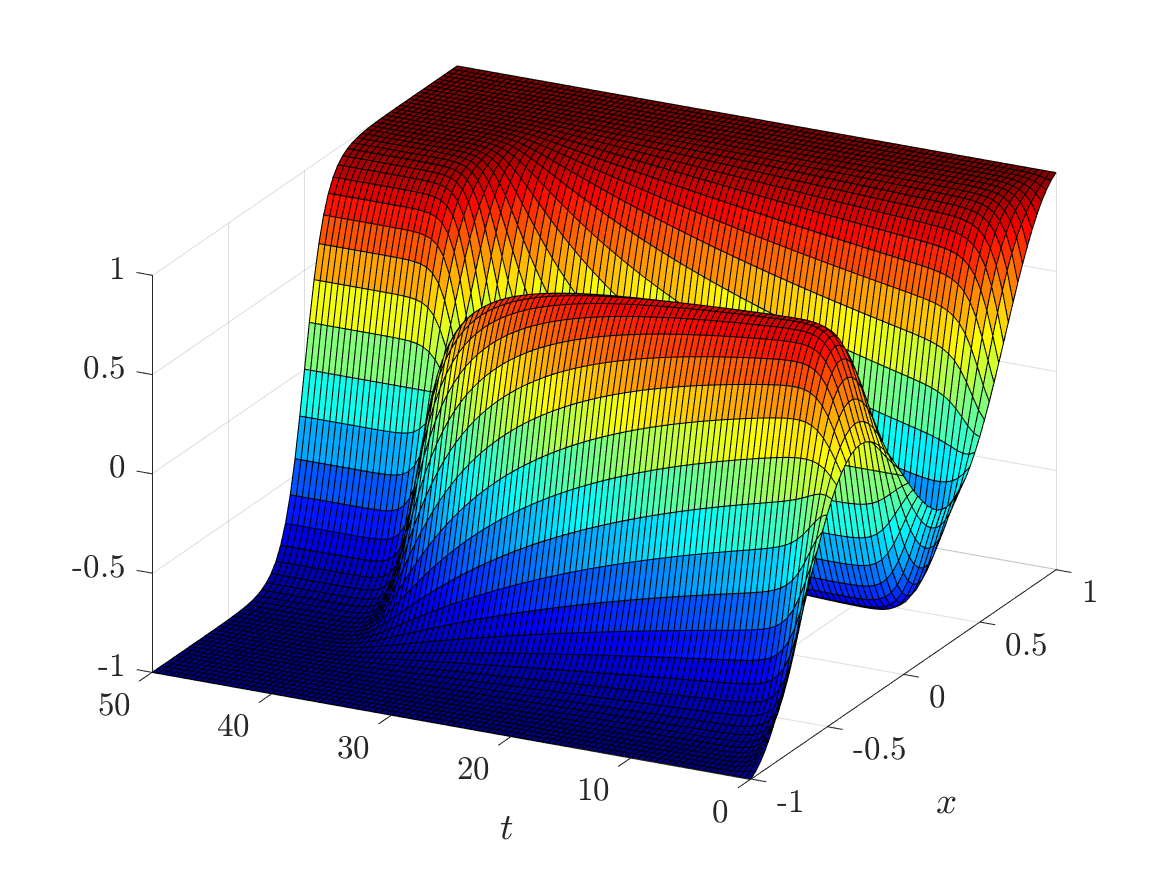}
}
\caption{Solution of the Allen-Cahn equation \eqref{AllenCahn} obtained with the three-stage DPG method (\ref{generalDPG3}).}
\label{AllenCahnSol}
\end{figure}

We consider the 1D+time Allen-Cahn equation \cite{feng2013nonlinear} in $\Omega=(-1,1)$
\begin{equation}\label{AllenCahn}
u_t=\epsilon\Delta u-f(u),
\end{equation}
with $f(u)=u^3-u$, $\epsilon=0.01$, $T=50$, non-homogeneous Dirichlet boundary conditions $u(-1,t)=-1$ and  $u(1,t)=1$, and initial condition $u_0(x)=0.53x+0.47\sin(-1.5\pi x)$. The Allen-Cahn equation \eqref{AllenCahn} can be obtained from the following Lyapunov energy functional
$$\mathcal{E}(u)=\int_{\Omega}\left(\frac{\epsilon}{2}|\nabla u|^2+\mathcal{F}(u)\right)dx,$$
where $\mathcal{F}(u)=\frac{1}{4}(u^2-1)^2$ is the Ginzburg-Landau double-well potential and $\mathcal{F}'(u)=f(u)$.
Figure \ref{AllenCahnSol} shows the solution of \eqref{AllenCahn} employing second-order central finite differences in space with $2^6$ points and the fourth-order DPG multistage method with $10^2$ time steps. Here, we can observe the \textit{metastability} phenomenon where the system exhibits a quasi-stable state that persists for a long time before transitioning to the true stable state. Finally, in Figure we conclude that the discrete energy decaying property, i.e., 
$$\mathcal{E}(u^n)\geq\mathcal{E}(u^{n+1}),\;n=0,\ldots,N-1,$$
is satisfied for all three methods.

\begin{figure}[h!]
\centering
\begin{tikzpicture}
\begin{axis}[small,
	width=1.75*\width,
	height=1.25*\height,
	xmin=0,
	xmax=50,
	xlabel=time,
	ylabel=Energy,
	grid=both,
    legend style={font=\small,at={(0.5,1.03)},anchor=south},
    legend entries={Hybrid Euler,DPG2,DPG3},
    legend columns=3
	]      
\addplot[
	color=black,dotted,line width=1.5pt, 
	]
table[x=t,y=Energy_EulerExpmv]{Results/EnergyAllenCahn_Ros_r6_m100.txt};
\addplot[
	color=blue,dotted,line width=1.5pt, 
	]
table[x=t,y=Energy_RK2Expmv]{Results/EnergyAllenCahn_Ros_r6_m100.txt};
\addplot[
	color=red,dotted,line width=1.5pt, 
	]
table[x=t,y=Energy_RK3Expmv]{Results/EnergyAllenCahn_Ros_r6_m100.txt};
\end{axis}
\end{tikzpicture}
\caption{Energy decaying property of the hybrid exponential Euler method (\ref{ExpEuler}), the two-stage DPG method (\ref{generalDPG2}), and the three-stage DPG method (\ref{generalDPG3}).}
\label{AllenCahnEnergy}
\end{figure}
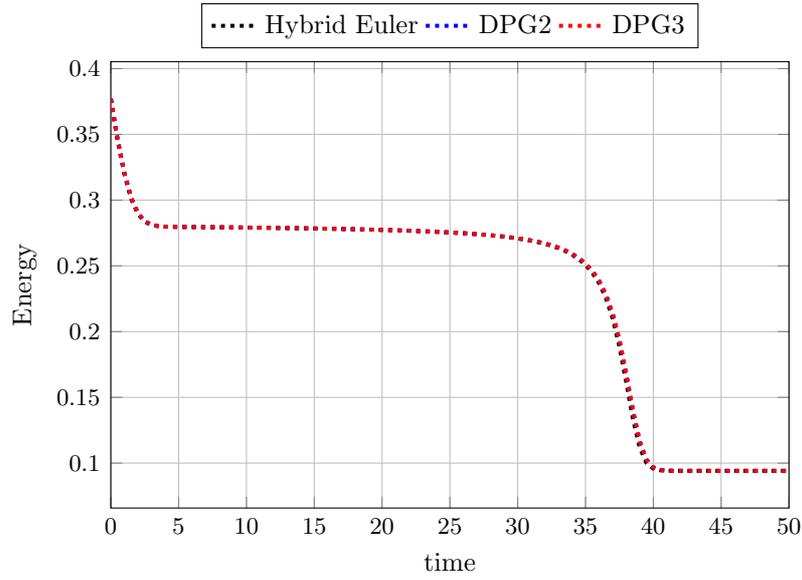

\subsection{Inviscid Burgers' equation}
\begin{figure}[h!]
\center{
\includegraphics[scale=0.5]{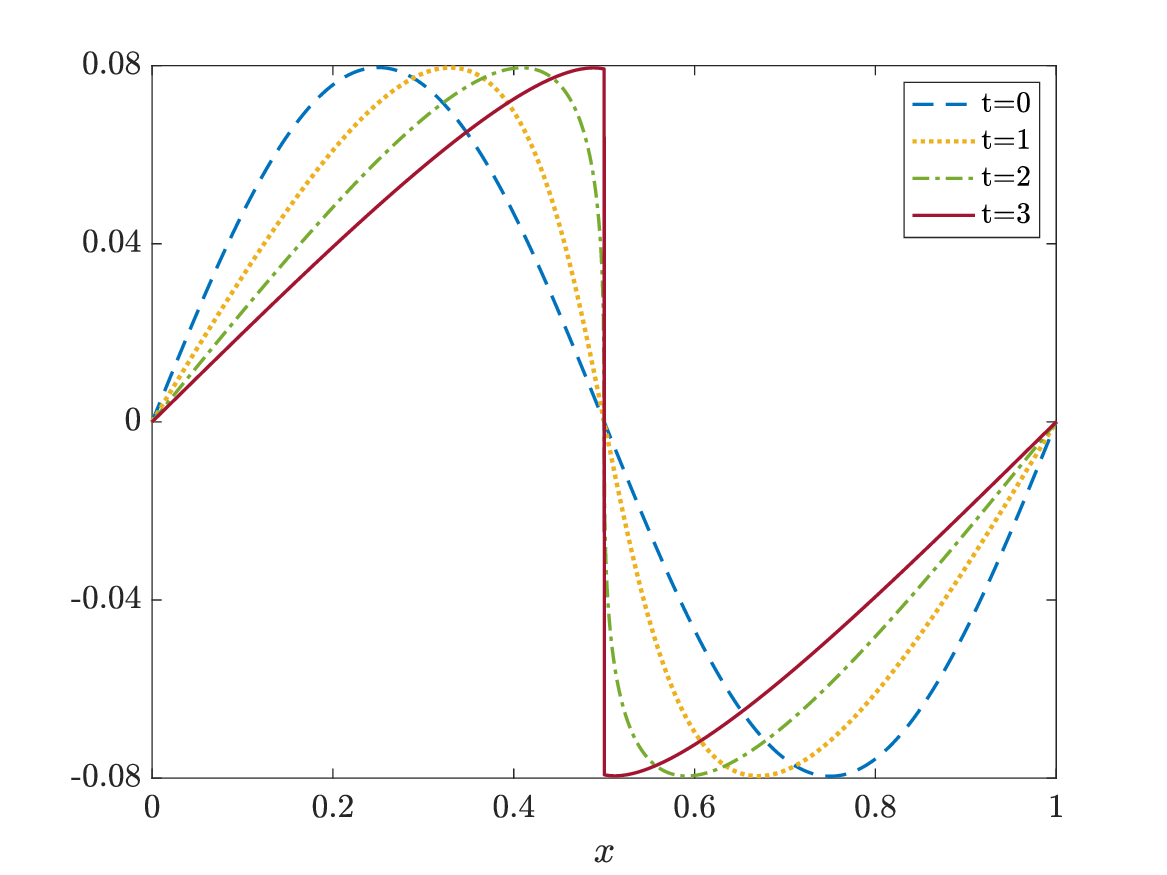}
}
\caption{Solution of the Burgers' equation \eqref{Burgers} at $t=0,1,2,3$ obtained with the three-stage DPG method (\ref{generalDPG3}).}
\label{BurgerSol}
\end{figure}

\begin{figure}[h!]
\centering
\begin{tikzpicture}
\begin{axis}[small,
	width=1.5*\width,
	height=1.25*\height,
	ymode=log,
	ymin=1e-3,
	ymax=3e-3,
	xlabel=t,
	ylabel=Kinetic energy,
	xmin=0,
	xmax=3,
	xtick={0,1,2,3},
	grid=both,
    legend style={font=\small,at={(0.5,1.03)},anchor=south},
    legend entries={Hybrid Euler,DPG2,DPG3},
    legend columns=3
	]      
\addplot[
	color=black,dotted,line width=1.5pt, 
	]
table[x=t,y=Energy_EulerExpmv]{Results/EnergyBurgers_Ros_r12_m10.txt};
\addplot[
	color=blue,dotted,line width=1.5pt,  
	]
table[x=t,y=Energy_RK2Expmv]{Results/EnergyBurgers_Ros_r12_m10.txt};
\addplot[
	color=red,dotted,line width=1.5pt,  
	]
table[x=t,y=Energy_RK3Expmv]{Results/EnergyBurgers_Ros_r12_m10.txt};
\end{axis}
\end{tikzpicture}
\caption{Kinetic energy at each time step for the hybrid exponential Euler method (\ref{ExpEuler}), the two-stage DPG method (\ref{generalDPG2}), and the three-stage DPG method (\ref{generalDPG3}).}
\label{EnergyBurgers}
\end{figure}
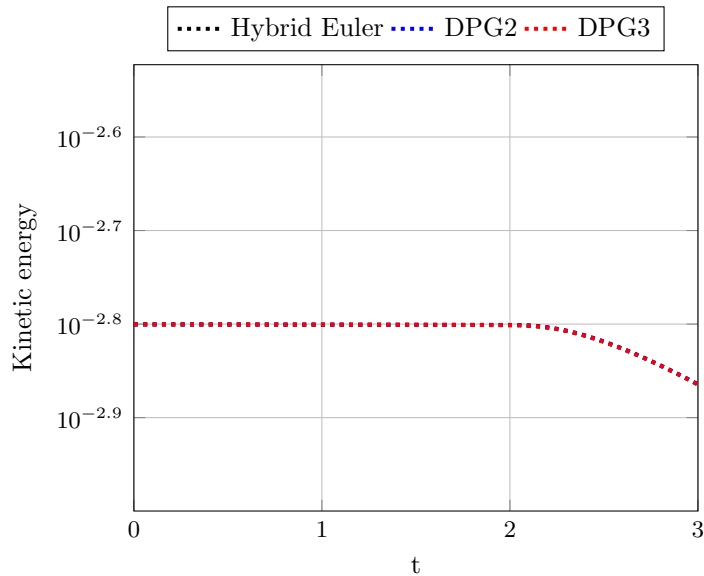

We consider the 1D+time in Inviscid Burgers' equation \cite{anguelov2008energy} in $\Omega=(0,1)$
\begin{equation}\label{Burgers}
u_t+uu_x=0,
\end{equation}
with $T=3$, periodic boundary conditions and initial condition $u_0(x)=\frac{1}{4\pi}\sin(2\pi x)$. Figure \ref{BurgerSol} shows several snapshots of the solution where we employed the Upwind scheme for space discretization and the fourth-order DPG multistage scheme. Here, we observe a wave propagating with characteristic speed $u$, and at $t=2$, there is a shock formation from the initially smooth solution. The Kinetic energy for equation \eqref{Burgers} can be defined as 
$$E(u(t))=\frac{1}{2}\lVert u(t)\rVert^2_0,$$
where $||\cdot ||_0$ denotes the $L^2$ error in $\Omega$. We know from \cite{anguelov2008energy} that before the formation of the shock, the solution is smooth, and it is satisfied that $E(u(t)) = E(u_0)$. Figure \ref{EnergyBurgers} shows that the methods we present in this article preserve the discrete Kinetic energy before the shock.

\section{Conclusions and future work}\label{Sec:Conclusions}
In this article, we derive three multistage methods up to order four employing the DPG time-marching scheme we developed in \cite{munoz2021equivalence,munoz2021adpg} for linear problems. We consider general nonlinear systems of ODEs and perform a linearization of the problem at each time step as in exponential Rosebrock methods. We first introduce the hybridization of the classical exponential Euler method, consisting of point values at each time step and piecewise constants inside the time interval. We observe that whereas we need an action of a $\varphi-$function to compute the interior, we only need a matrix-vector product to compute the trace variable. We then employ this hybrid Euler method for the internal stages to construct the third and fourth-order methods and derive the stiff-order conditions to determine the values of the coefficient functions. We observe that for the fourth-order method, we need to additionally introduce a correction term that depends upon the Jacobian. Finally, we provide a convergence proof under classical stability conditions on the exponential of the Jacobian matrix that changes at each time step. For the convergence proof, we ensure that the correction factor is bounded by the error at the internal stages. We test the numerically the convergence and performance predicted by the theory on the so-called 2D+time Hochbruck-Ostermann equation. Finally, we test the energy decaying property on the 1D+time Allen-Cahn equation and the Kinetic energy preservation property on the 1D+time Burgers' equation. We conclude that the proposed linearization does not break the intrinsic structure of the quasilinear models considered in this work.

In the literature, the highest-order exponential Rosenbrock method is of order five \cite{luan2014exponential}. We are interested in seeing if we can derive multistage methods employing the DPG construction for the internal stages to go beyond this threshold. Therefore, we will study in future work the use of higher-order polynomials, for example, piecewise linears or quadratics as we introduced in \cite{munoz2021equivalence}, to derive multistage DPG methods of order higher than four.

\appendix 
\section{Summary of notation}\label{App:Summarize}
The following tables summarize the meaning of each symbol and errors employed throughout the convergence proofs in Sections \ref{Sec:Stiff} and \ref{Sec:Convergence}. 
\renewcommand{\arraystretch}{1.5}
\begin{table}[h!]
\centering
\begin{tabular}{|c|c|c|c|}\hline
$u_n$& numerical approximation of the solution at $t_n$\\ \hline
$u_{n_2}$, $u_{n_3}$& internal stages computed with $u_n$\\ \hline
$u_{n+1}$& update of the numerical scheme computed with $u_n$, $u_{n_2}$ and $u_{n_3}$\\ \hline
$J_n$& Jacobian evaluated in $u_n$\\ \hline
$g_n$& nonlinear remainder after linearization at $u_n$\\ \hline
$C_{n_2}$& correction factor computed with $g'_n$, $u_n$, $u_{n_2}$ and $u_{n_3}$\\ \hline
$\hat{u}_n$, $\hat{u}_{n+1}$& analytical solution at $t_n$ and $t_{n+1}$\\ \hline
$\bar{u}_{n_2}$, $\bar{u}_{n_3}$& internal stages computed with $\hat{u}_n$\\ \hline
$\bar{u}_{n+1}$& update of the numerical scheme computed with $\hat{u}_n$, $\bar{u}_{n_2}$ and $\bar{u}_{n_3}$\\ \hline
$\hat{J}_n$& Jacobian evaluated in $\hat{u}_n$\\ \hline
$\hat{g}_n$& nonlinear remainder after linearization at $\hat{u}_n$\\ \hline
$\bar{C}_{n_2}$& correction factor computed with $\hat{g}'_n$, $\hat{u}_n$, $\bar{u}_{n_2}$ and $\bar{u}_{n_3}$\\ \hline
\end{tabular}
\caption{Summary of symbols.}
\label{SummarySymbols}
\end{table}
\begin{table}[h!]
\centering
\begin{tabular}{|c|c|c|c|}\hline
$e_n=u_n-\hat{u}_n$\\ \hline
$\bar{e}_{n+1}=u_{n+1}-\bar{u}_{n+1}$\\ \hline
$\hat{e}_{n+1}=\bar{u}_{n+1}-\hat{u}_{n+1}$\\ \hline
$e_{n+1}=\bar{e}_{n+1}+\hat{e}_{n+1}=u_{n+1}-\hat{u}_{n+1}$\\ \hline
$E_{n_2}=u_{n_2}-\bar{u}_{n_2}$, $E_{n_3}=u_{n_3}-\bar{u}_{n_3}$\\ \hline
\end{tabular}
\caption{Summary of errors.}
\label{SummaryErrors}
\end{table}

\section*{Acknowledgements}
Judit Mu\~noz-Matute has received funding from the European Union's Horizon 2020 research and innovation programme under the Marie Sklodowska-Curie individual fellowship grant agreement No. 101017984 (GEODPG). Leszek  Demkowicz was partially supported with NSF grant No. 1819101.

\section*{References}
\bibliographystyle{abbrv} 
\bibliography{references.bib}

\end{document}